\documentclass[11pt]{article}
\usepackage{amssymb,amsmath,booktabs,hyperref,mathtools}
\usepackage{graphics,authblk}
\usepackage{xcolor}
\usepackage{longtable,nicefrac}
\usepackage{graphicx}
\usepackage{subcaption}
\usepackage{multirow}
\usepackage{tabularx,enumerate}
\usepackage{rotating}
\usepackage{dcolumn,accents}
\usepackage{pdflscape}
\usepackage{array}
\usepackage{lscape}
\usepackage{anysize}
\usepackage{color}
\usepackage{amsthm}
\usepackage{listings}
\newcommand{\mathleft}{\@fleqntrue\@mathmargin0pt}
\newcommand{\mathcenter}{\@fleqnfalse}

\providecommand{\keywords}[1]{\textbf{\textit{Keywords: }} #1}
\usepackage{float}
\usepackage{soulutf8}
\usepackage{adjustbox}
\usepackage{colortbl}
\usepackage{blindtext,nameref}

\newtheorem{theorem}{Theorem}

\newtheorem*{ivt}{Intermediate Value Theorem (IVT)}
\newtheorem{proposition}{Proposition}

\theoremstyle{definition}
\newtheorem{definition}{Definition}
\theoremstyle{remark}

\graphicspath{{photos/}}
\allowdisplaybreaks
\begin{document}
	\title{A Generalized Analytical Framework for the Nonlinear Best-Worst Method}
	\author{Harshit M. Ratandhara, Mohit Kumar}
	\date{}
	\affil{Department of Basic Sciences,\\ Institute of Infrastructure, Technology, Research And Management, Ahmedabad, Gujarat-380026, India\\ Email: harshitratandhara1999@gmail.com, mshramadma.iitr@gmail.com}
	\maketitle
	\begin{abstract}
		The nonlinear model of the best-worst method frequently produces multiple optimal weight sets, which are conventionally determined through optimization software. While an analytical approach exists that provides both a closed-form expression for the optimal interval-weights and a secondary objective function to determine the best optimal weight set, we demonstrate that this approach is only valid when preferences are quantified using the Saaty scale and only a single decision-maker is involved. To tackle this issue, we propose a framework compatible with any scale and any number of decision-makers. We first derive an analytical expression for optimal interval-weights and then select the best optimal weight set. After demonstrating that the values of consistency index for the Saaty scale in the existing literature are not well-defined, we derive a formula of consistency index. We also obtain an analytical expression for the consistency ratio, enabling its use as an input-based consistency indicator. Furthermore, we establish that when multiple best/worst criteria are present, weights may vary among best criteria and among the worst criteria. To address this limitation, we modify the original optimization model for weight computation in such instances, solve it analytically to obtain optimal interval-weights and then select the best optimal weight set using a secondary objective function. Finally, we demonstrate and validate the proposed approach using numerical examples and a real-world case study of ranking barriers to energy efficiency in buildings.
	\end{abstract}
	\section{Introduction}
	Multi-Criteria Decision-Making (MCDM), a fundamental branch of operations research, addresses complex decision scenarios involving conflicting criteria. Typically, MCDM approaches solve decision problems through two sequential steps: determining weights of decision criteria and ranking of alternatives. This classification divides MCDM methods into two categories: weight determination methods (such as AHP \cite{saaty1994make}, ANP \cite{saaty2004decision}, and SMARTS \cite{edwards1994smarts}) and ranking methods (including TOPSIS \cite{hwang1981methods}, VIKOR \cite{opricovic2004compromise}, ELECTRE \cite{roy1991outranking}, and PROMETHEE \cite{brans1985note}).\\\\
	The Best-Worst Method (BWM) is a widely used weight calculation method that utilizes pairwise comparisons between decision criteria to derive weights \cite{rezaei2015best}. In this method, the optimal weights are obtained by minimizing the distance between weight ratios and given comparison values. Based on different distance functions, various models of BWM have been developed. The original model of BWM, proposed by Rezaei \cite{rezaei2015best}, uses maximum deviation as the distance function. Since this approach involves solving a nonlinear optimization problem, it is known as the nonlinear BWM. Kocak et al. \cite{kocak2018euclidean} introduced a model of BWM based on Euclidean distance, termed the Euclidean BWM. Brunelli and Rezaei \cite{brunelli2019multiplicative} proposed the multiplicative BWM, incorporating the metric $\max\{\nicefrac{x}{y},\nicefrac{y}{x}\}$ defined on $(\mathbb{R_+},\cdot,\leq)$ into the BWM framework. Amiri and Emamat \cite{amiri2020goal} developed a goal programming-based BWM using the taxicab distance. Tu et al. \cite{tu2023priority} introduced two alternative BWM formulations: the approximate eigenvalue model and the logarithmic least squares model. For group decision-making, Safarzadeh et al. \cite{safarzadeh2018group} proposed two models of BWM, one based on total deviation and the other based on maximum deviation. Xu and Wang \cite{xu2024some} discussed eleven BWM models for individual Decision-Makers (DMs) and nine for group decision-making contexts. Mohammadi and Rezaei \cite{mohammadi2020bayesian} incorporated Bayesian into the BWM framework, establishing a probabilistic approach to group decision-making. Corrente et al. \cite{corrente2024better} developed the parsimonious BWM, an enhanced version of the nonlinear BWM specifically designed for decision contexts involving numerous alternatives.\\\\
	To address the non-uniqueness of optimal weight sets in the nonlinear BWM, Rezaei \cite{rezaei2016best} derived optimal interval-weights using two optimization models. Later, he introduced the concentration ratio to measure the dispersion of these interval-weigths \cite{rezaei2020concentration}. He also developed the linear BWM \cite{rezaei2016best}, which retains the underlying philosophy of the nonlinear BWM while transforming its optimization framework into a linear formulation. While this model guarantees a unique weight set, its feasible region differs from the original nonlinear approach. Wu et al. \cite{wu2023analytical} developed an analytical approach to derive optimal interval-weights for the nonlinear BWM, eliminating the model's dependency on optimization software. Building on this foundation, they introduced a secondary objective function to determine the best optimal weight set from the solution space. Ratandhara and Kumar \cite{ratandhara2024analytical} subsequently proposed an analytical framework for the multiplicative BWM, achieving the same objectives of software independence and selection of the best optimal weight set.\\\\
	Measurement of consistency of decision data is crucial in MCDM methods since outcomes depend directly on this input data. In BWM, consistency evaluation is performed through the Consistency Ratio (CR), which is computed using the optimal objective value of optimization model and Consistency Index (CI) \cite{rezaei2015best}. This output-based consistency indicator can only provide feedback about inconsistencies after completing all calculations, resulting in reduced time efficiency. Liang et al. \cite{liang2020consistency} introduced an alternative input-based consistency indicator called input-based CR, along with establishing threshold values for both output-based and input-based CR to check admissibility of preference values. Furthermore, Lei et al. \cite{lei2022preference} developed an optimization model that recommends optimal preference modifications while achieving both ordinal consistency and an acceptable level of cardinal consistency.\\\\
	To handle uncertain preferences, several fuzzy-set-based extensions of BWM have been proposed. Guo and Zhao \cite{guo2017fuzzy} extended the nonlinear BWM to a fuzzy environment, while Rostami et al. \cite{rostami2023goal} introduced a fuzzy adaptation of the goal programming-based BWM. Additionally, Ratandhara and Kumar \cite{ratandhara2024alpha} proposed an $\alpha$-cut interval-based model of fuzzy BWM. The BWM framework has also been extended to more advanced uncertainty representations, including intuitionistic fuzzy sets \cite{wan2021novel,wan2024novel,cheng2024decision}, hesitant fuzzy sets \cite{ali2019hesitant,li2021approaches} and spherical fuzzy sets \cite{haseli2024extension}. Moreover, BWM has been integrated with other MCDM techniques, such as BWM-VIKOR \cite{aleksic2023evaluation,gao2024novel}, BWM-ELECTRE \cite{yadav2018hybrid,chen2024selecting}, BWM-TOPSIS \cite{varchandi2024integrated,trivedi2024hybrid}, BWM-MULTIMOORA \cite{hafezalkotob2019interval,yucesan2024evaluating} and Best-Worst Tradeoff (BWT) method \cite{liang2022best}. Owing to its simplicity and reliability, BWM has been widely applied in practical decision-making problems, notably in location selection \cite{uyan2023gis,amjad2024identification}, logistics risk assessment \cite{delice2020new,ding2023assessment}, and supplier selection \cite{tavana2024integrated,varchandi2024integrated}.\\\\
	In this study, we establish the following key research gaps in the nonlinear BWM.
	\begin{enumerate}[(i)]
		\item The framework proposed by Wu et al. \cite{wu2023analytical} works well when preferences are quantified using the Saaty scale, but leads to different optimal interval-weights and a best optimal weight set than their actual values (in fact, non-well-defined optimal interval-weights) when some other scale such as the Salo-H{\"a}m{\"a}l{\"a}inen, Lootsma or Donegan-Dodd-McMaster scale is used. This is a significant concern, as the quantification of preferences using the Saaty scale may not always be reliable \cite{koczkodaj2016pairwise}.
		\item The existing framework leads to non-well-defined optimal interval-weights and a different best weight set than the actual one in the presence of multiple DMs even when preferences are quantified using the Saaty scale.
		\item The values of CI for the Saaty scale computed by Rezaei \cite{rezaei2015best} are not well-defined.
		\item In the presence of multiple best/worst criteria, computing weights by arbitrarily selecting any one best/worst criterion leads to weight differences (in both interval-weights and the best weights) among the best criteria and among the worst criteria.
	\end{enumerate}
	In this work, we propose an analytical framework for the nonlinear BWM that is compatible with any preference scale and any number of DMs. Our approach derives optimal interval-weights by solving an optimal modification-based optimization problem, which is equivalent to the original optimization model. From the collection of all optimal weight sets, we also select the best optimal weight set by introducing a secondary objective function. Based on this framework, we derive a formula for computing CI that is valid for any scale, along with an analytical expression for CR that serves as an input-based consistency indicator. We further establish essential properties of CR to verify its reliability as a consistency measure. To ensure equal weights for all best and all worst criteria in cases with multiple best/worst criteria, we develop a modified optimization model and analytically derive both the corresponding optimal interval-weights and the best optimal weight set. We validate the proposed approach through numerical examples. Additionally, we rank barriers to energy efficiency in buildings to demonstrate applicability of the framework.\\\\
	The remainder of this paper is organized as follows. Section 2 introduces fundamental definitions and provides a concise overview of the nonlinear BWM and its existing framework. Section 3 identifies and examines critical research gaps in the current methodology. Section 4 presents a generalized analytical framework for the nonlinear BWM including weight calculation, consistency analysis, and numerical examples. Section 5 demonstrates the applicability of the proposed approach by ranking barriers to energy efficiency in buildings. Finally, Section 6 discusses concluding remarks and suggests potential directions for future research.
	\section{Preliminaries}
	In this section, we discuss some fundamental definitions and results, along with a brief overview of the nonlinear BWM and its analytical framework.
	\subsection{Basic Concepts and Results}
	Let $C=\{c_1,c_2,\ldots,c_n\}$ be the set of decision criteria, and let $D=\{c_1,c_2,\ldots,c_n\}\setminus\{c_{b},c_{w}\}$ throughout the article. For brevity, when the context is clear, we simplify the notation by representing the sets as $C=\{1,2,\ldots,n\}$ and $D=\{1,2,\ldots,n\}\setminus\{b,w\}$.\\\\
	The Pairwise Comparison System (PCS) is the pair $(A_b,A_w)$, where $A_b=(a_{b1},a_{b2},\ldots,a_{bn})$ and $A_w=(a_{1w},a_{2w},\ldots,a_{nw})^T$ denote the best-to-other and the other-to-worst vector respectively. Here, $a_{ij}$ represents the relative preference of the $i^{th}$ criterion over the $j^{th}$ criterion.
	\begin{definition} \cite{rezaei2015best}
		A PCS $(A_b,A_w)$ is said to be consistent if $a_{bi}\times a_{iw}=a_{bw}$ for all $i\in D$.
	\end{definition}
	\begin{theorem}\cite{wu2023analytical}\label{5accurate_weight}
		The system of linear equations
		\begin{equation}\label{5system}
			\begin{split}
				&\frac{w_b}{w_i}=a_{bi},\quad \frac{w_i}{w_w}=a_{iw},\quad \frac{w_b}{w_w}=a_{bw},\ i\in D,\\
				&w_1+w_2+\ldots+w_n=1
			\end{split}
		\end{equation}
		has a solution if and only if $(A_b,A_w)$ is consistent. Also, if solution exists, then it is unique and is given by
		\begin{equation}\label{accurate_weights}
			w_i=\frac{a_{iw}}{\displaystyle\sum_{j\in C}a_{jw}}=\frac{1}{a_{bi}\displaystyle\sum_{j\in C}\frac{1}{a_{bj}}},\ i\in C.
		\end{equation}
	\end{theorem}
	\begin{ivt}\cite{royden1988real}
		Let $f$ be a real-valued continuous function on $[a,b]$ for which $f(a)<c<f(b)$. Then there is a $a<x_0<b$ such that $f(x_0)=c$.
	\end{ivt}
	\subsection{Nonlinear BWM}
	BWM is a pairwise comparison-based MCDM method for deriving criteria weights \cite{rezaei2015best}. The steps of BWM are as follows.\\\\
	\textbf{Step 1:} Formation of the set of decision criteria $C$.\\\\
	\textbf{Step 2:} Selection of the best criterion $c_b$ and the worst criterion $c_w$ from $C$.\\\\
	\textbf{Step 3:} Determination of the PCS $(A_b,A_w)$.\\\\
	The preferences $a_{ij}$ are generally expressed as linguistic terms, which are then quantified using an established scale. Some such scales are given in Table \ref{5scale}.
	\begin{table}[H]
		\caption {Quantification of linguistic terms using different scales \label{5scale}}
		\centering
		\begin{tabular}{@{}ccccc@{}}
			\toprule[0.1em]
			\multirow{2}{*}{Linguistic term}&Saaty&Salo-H{\"a}m{\"a}l{\"a}inen&Lootsma&Donegan-Dodd-McMaster\\
			&scale \cite{saaty1994make}&scale \cite{ji2003scale}&scale \cite{ji2003scale}&scale (7-based) \cite{dodd1995scale}\\
			\midrule
			Indifference&$1$&$1$&$1$&$1$\\
			-&$2$&$1.2222$&$\sqrt{2}$&$1.1257$\\
			Moderate preference&$3$&$1.5$&$2$&$1.2715$\\
			-&$4$&$1.8571$&$2\sqrt{2}$&$1.4470$\\
			Strong preference&$5$&$2.3333$&$4$&$1.6684$\\
			-&$6$&$3$&$4\sqrt{2}$&$1.9670$\\
			Very strong preference&$7$&$4$&$8$&$2.4142$\\
			-&$8$&$5.6667$&$8\sqrt{2}$&$3.2289$\\
			Extreme preference&$9$&$9$&$16$&$5.8284$\\
			\bottomrule[0.1em]
		\end{tabular}
	\end{table}
	\hspace{-0.7cm}
	\textbf{Step 4:} Computation of optimal weights using a nonlinear optimization model.\\
	Consider the following minimization problem.
	\begin{equation}\label{5optimization_1}
		\begin{split}
			&\min\epsilon \\
			&\text{subject to:}\\	
			&\left|\frac{w_b}{w_i}-a_{bi}  \right| \leq \epsilon, \quad \left|\frac{w_i}{w_w}-a_{iw}  \right| \leq \epsilon, \quad	\left| \frac{w_b}{w_w}-a_{bw} \right| \leq \epsilon, \\
			&w_1+w_2+\ldots+w_n=1,\quad w_j\geq 0\text{ for all } i\in D \text{ and } j\in C.
		\end{split}
	\end{equation}
	Problem \eqref{5optimization_1} has optimal solution(s) of the form $(w_1^*,w_2^*,\ldots,w_n^*,\epsilon^*)$. For each optimal solution, $W^*=\{w_1^*,w_2^*,\ldots,w_n^*\}$ is an optimal weight set, while $\epsilon^*$ indicates its accuracy. Since $\epsilon^*$ is the optimal objective value, it remains the same for all optimal weight sets.\\\\
	To address the non-uniqueness of optimal solutions in problem \eqref{5optimization_1}, Rezaei \cite{rezaei2016best} employed interval-analysis, observing that the set of all optimal weights for each criterion is an interval. These optimal interval-weights can be obtained using the following optimization problems.
	\begin{equation}\label{5optimization_3}
		\begin{split}
			&\min w_k \\
			&\text{subject to:}\\	
			&\left|\frac{w_b}{w_i}-a_{bi}  \right| \leq \epsilon^*, \quad \left|\frac{w_i}{w_w}-a_{iw}  \right| \leq \epsilon^*, \quad	\left| \frac{w_b}{w_w}-a_{bw} \right| \leq \epsilon^*, \\
			&w_1+w_2+\ldots+w_n=1,\quad w_j\geq 0\text{ for all } i\in D \text{ and } j\in C.
		\end{split}
	\end{equation}
	\begin{equation}\label{5optimization_4}
		\begin{split}
			&\max w_k \\
			&\text{subject to:}\\	
			&\left|\frac{w_b}{w_i}-a_{bi}  \right| \leq \epsilon^*, \quad \left|\frac{w_i}{w_w}-a_{iw}  \right| \leq \epsilon^*, \quad	\left| \frac{w_b}{w_w}-a_{bw} \right| \leq \epsilon^*, \\
			&w_1+w_2+\ldots+w_n=1,\quad w_j\geq 0\text{ for all } i\in D \text{ and } j\in C.
		\end{split}
	\end{equation}
	Problems \eqref{5optimization_3} and \eqref{5optimization_4} are optimization problems having $n$ variables, where the Greatest Lower Bound (GLB) and the Least Upper Bound (LUB) of the optimal interval-weight for criterion $c_k$ serve as the respective optimal objective values, i.e., if $w_k'^*$ and $w_k''^*$ denote the optimal objective values of problems \eqref{5optimization_3} and \eqref{5optimization_4} respectively, then the optimal interval-weight for $c_k$ is $[w_k'^*,w_k''^*]$.\\\\
	The effectiveness of an MCDM method depends on the decision data that is often inconsistent because of human engagement. A key requirement for any MCDM method is the ability to assess and quantify these inconsistencies. In BWM, this assessment is performed using the Consistency Ratio (CR) defined as
	\begin{equation}\label{5CR}
		\text{CR}=\frac{\epsilon^*}{\text{Consistency Index (CI)}},
	\end{equation}
	where CI $=\sup\{\epsilon^*:\epsilon^*$ is the optimal objective value of problem \eqref{5optimization_1} for some $(A_b,A_w)$ having the given value of $a_{BW}$\} \cite{rezaei2015best}. So, CI is a function of $a_{bw}$. The values of CI for the Saaty scale are given in Table \ref{5ci}.
	\begin{table}[H]
		\caption {The values of CI for the Saaty scale \cite{rezaei2015best}\label{5ci}}
		\centering
		\begin{tabular}{@{}ccccccccc@{}}
			\toprule[0.1em]
			$a_{bw}$&$2$&$3$&$4$&$5$&$6$&$7$&$8$&$9$\\
			\midrule
			CI&$0.4384$&$1$&$1.6277$&$2.2984$&$3$&$3.7250$&$4.4688$&$5.2279$\\
			\bottomrule[0.1em]
		\end{tabular}
	\end{table}
	\subsection{An Analytical Framework for the Nonlinear BWM}
	Wu et al. \cite{wu2023analytical} proposed an analytical approach to derive optimal interval-weights without requiring optimization software. They also introduced a secondary objective function to select the best optimal weight set from the collection of all optimal weight sets.
	\subsubsection{Interval-Weights}
	The optimal interval-weights are calculated using the following optimization model, which is driven by an optimal modification of the PCS.
	\begin{equation}\label{5optimization_2}
		\begin{split}
			&\min\eta \\
			&\text{subject to:}\\	
			&\left|\tilde{a}_{bi}-a_{bi}  \right| \leq \eta, \quad \left|\tilde{a}_{iw}-a_{iw}  \right| \leq \eta, \quad	\left| \tilde{a}_{bw}-a_{bw} \right| \leq \eta, \\
			&\tilde{a}_{bi}\times\tilde{a}_{iw}=\tilde{a}_{bw},\quad \tilde{a}_{bi},\tilde{a}_{iw},\tilde{a}_{bw}\geq 0\text{ for all } i\in D.
		\end{split}
	\end{equation}
	Problem \eqref{5optimization_2} has optimal solution(s) of the form $(\tilde{a}_{bi}^*,\tilde{a}_{iw}^*,\tilde{a}_{bw}^*,\eta^*)$, where $i\in D$. Each optimal solution, along with $\tilde{a}_{bb}^*=\tilde{a}_{ww}^*=1$, leads to a consistent PCS $(\tilde{A}_b^*,\tilde{A}_w^*)$, referred to as the optimally modified PCS and $\eta^*$ indicates its the accuracy. Since $\eta^*$ is the optimal objective value, it remains the same for all $(\tilde{A}_b^*,\tilde{A}_w^*)$.\\\\
	Wu et al. \cite{wu2023analytical} established that
	\begin{equation}\label{5equivalence}
		\begin{minipage}{\textwidth}
			\begin{enumerate}
				\item  $\epsilon^*=\eta^*$
				\item for each $W^*=\{w_1^*,w_2^*,\ldots,w_n^*\}$, there exists a unique $(\tilde{A}_b^*,\tilde{A}_w^*)$ satisfying the relation $w_i^*=\frac{\tilde{a}_{iw}^*}{\sum_{i=1}^{n}\tilde{a}_{jw}^*}$.
			\end{enumerate}
		\end{minipage}
	\end{equation}
	This implies that problem \eqref{5optimization_2} is equivalent to problem \eqref{5optimization_1}. Consequently, the analytical expression for optimal solutions of problem \eqref{5optimization_1} can be obtained by solving problem \eqref{5optimization_2} analytically. To describe the analytical solution of problem \eqref{5optimization_2}, some mathematical symbols must first be introduced.\\\\
	Let $D_1=\{i\in D: a_{bi}\times a_{iw}<a_{bw}\}$, $D_2=\{i\in D: a_{bi}\times a_{iw}>a_{bw}\}$, and $D_3=\{i\in D: a_{bi}\times a_{iw}=a_{bw}\}$.\\\\
	Fix $i\in D$. Then there are three possibilities.
	\begin{enumerate}[(i)]
		\item $i \in D_1$\\\\
		Consider the quadratic equation
		\begin{equation}\label{5quad_1}
			(a_{bi}+x)\times (a_{iw}+x)=a_{bw}-x.
		\end{equation}
		Let $f(x)=(a_{bi}+x)\times (a_{iw}+x)$ and $g(x)=a_{bw}-x$, where $x\in \mathbb{R}$. Note that $f(0)=a_{bi}\times a_{iw}< a_{bw} =g(0)$ and $f(a_{bw})=(a_{bi}+a_{bw})\times (a_{iw}+a_{bw})>0 = g(a_{bw})$. So, by IVT, there exist $0<c<a_{bw}$ such that $f(c)=g(c)$, i.e., $c$ is a positive root of equation \eqref{5quad_1}. Let $\epsilon_i$ be the smallest positive root of equation \eqref{5quad_1}. Then
		\begin{equation}\label{5CV1}
			(a_{bi}+\epsilon_i)\times(a_{iw}+\epsilon_i)=a_{bw}-\epsilon_i.
		\end{equation}
		From the above discussion, it follows that $\epsilon_i<a_{bw}$.
		\item $i\in D_2$\\\\
		Consider the quadratic equation
		\begin{equation}\label{5quad_2}
			(a_{bi}-x)\times (a_{iw}-x)=a_{bw}+x.
		\end{equation}
		Let $f(x)=(a_{bi}-x)\times (a_{iw}-x)$ and $g(x)=a_{bw}+x$, where $x\in \mathbb{R}$. Let $a=\min\{a_{bi},a_{iw}\}$ Note that $f(0)=a_{bi}\times a_{iw}> a_{bw} =g(0)$ and $f(a)=0< a_{bw}+a= g(a)$. So, by IVT, there exist $0<c<a$ such that $f(c)=g(c)$, i.e., $c$ is a positive root of equation \eqref{5quad_2}. Let $\epsilon_i$ be the smallest positive root of equation \eqref{5quad_2}. Then
		\begin{equation}\label{5CV2}
			(a_{bi}-\epsilon_i)\times(a_{iw}-\epsilon_i)=a_{bw}+\epsilon_i.
		\end{equation}
		From the above discussion, it follows that $\epsilon_i<a$, i.e., $\epsilon_i<a_{bi}$ and $\epsilon_i<a_{iw}$.
		\item $i\in D_3$\\\\ 
		In this case, take $\epsilon_i=0$.
	\end{enumerate}
	So, in any case, we get
	\begin{equation}\label{5CV}
		\epsilon_i=\bigg|\frac{a_{bi}+a_{iw}+1-\sqrt{(a_{bi}+a_{iw}+1)^2-4(a_{bi}\times a_{iw}-a_{bw})}}{2}\bigg|.
	\end{equation}
	Now, fix $i,j\in D$. Then there are three possibilities.
	\begin{enumerate}[(i)]
		\item $a_{bi}\times a_{iw}<a_{bj}\times a_{jw}$\\\\
		In this case, take $\epsilon_{i,j}=\frac{a_{bj}\times a_{jw}-a_{bi}\times a_{iw}}{a_{bi}+a_{iw}+a_{bj}+a_{jw}}.$ This gives
		\begin{equation}\label{5CV3}
			(a_{bi}+\epsilon_{i,j})\times (a_{iw}+\epsilon_{i,j})= (a_{bj}-\epsilon_{i,j})\times(a_{jw}-\epsilon_{i,j}).
		\end{equation}
		Note that $\epsilon_{i,j}<a_{bj}$ and $\epsilon_{i,j}<a_{jw}$.
		\item $a_{bi}\times a_{iw}>a_{bj}\times a_{jw}$\\\\
		In this case, take $\epsilon_{i,j}=\frac{a_{bi}\times a_{iw}-a_{bj}\times a_{jw}}{a_{bi}+a_{iw}+a_{bj}+a_{jw}}.$ This gives
		\begin{equation}\label{5CV4}
			(a_{bi}-\epsilon_{i,j})\times (a_{iw}-\epsilon_{i,j})= (a_{bj}+\epsilon_{i,j})\times(a_{jw}+\epsilon_{i,j}).
		\end{equation}
		Note that $\epsilon_{i,j}<a_{bi}$ and $\epsilon_{i,j}<a_{iw}$.
		\item $a_{bi}\times a_{iw}=a_{bj}\times a_{jw}$\\\\
		In this case, take $\epsilon_{i,j}=0$.
	\end{enumerate}
	So, in any case, we get
	\begin{equation}\label{5CV'}
		\epsilon_{i,j}=\bigg|\frac{a_{bi}\times a_{iw}-a_{bj}\times a_{jw}}{a_{bi}+a_{iw}+a_{bj}+a_{jw}}\bigg|.
	\end{equation}
	Let $i_1\in D_1$ and $i_2\in D_2$ be such that $\epsilon_{i_1}=\max\{\epsilon_i:i\in D_1\}$ and $\epsilon_{i_2}=\max\{\epsilon_i:i\in D_2\}$. Then, by \cite[Proposition 3]{wu2023analytical}, we get
	\begin{equation}\label{5op_obj}
		\epsilon^*=\begin{cases}
			\epsilon_{i_1}\quad \ \ \ \text{if } (a_{bi_2}-\epsilon_{i_1})\times(a_{i_2w}-\epsilon_{i_1})\leq a_{bw}-\epsilon_{i_1},\\
			\epsilon_{i_2}\quad \ \ \ \text{if } (a_{bi_1}+\epsilon_{i_2})\times(a_{i_1w}+\epsilon_{i_2})\geq a_{bw}+\epsilon_{i_2},\\
			\epsilon_{i_1,i_2}\quad\text{otherwise}.
		\end{cases}
	\end{equation}
	Now, by \cite[Theorem 3]{wu2023analytical}, the collection of all optimally modified PCS is
	\begin{subequations}\label{5op_PCS}
		\begin{align}
			&\tilde{a}_{bw}^*=\begin{cases}
				a_{bw}-\epsilon^*\quad\quad\quad\quad\quad\quad\quad\quad\quad\quad\text{if } \epsilon^*=\epsilon_{i_1},\\
				a_{bw}+\epsilon^*\quad\quad\quad\quad\quad\quad\quad\quad\quad\quad\text{if } \epsilon^*=\epsilon_{i_2},\\
				(a_{bi_2}-\epsilon_{i_1,i_2})\times(a_{i_2w}-\epsilon_{i_1,i_2})\quad\text{if } \epsilon^*=\epsilon_{i_1,i_2},
			\end{cases} \label{5op_PCS_1}\\
			&\tilde{a}_{iw}^*\in \bigg[\max\biggl\{a_{iw}-\epsilon^*,\frac{\tilde{a}_{bw}^*}{a_{bi}+\epsilon^*}\biggr\},\min\biggl\{a_{iw}+\epsilon^*,\frac{\tilde{a}_{bw}^*}{a_{bi}-\epsilon^*}\biggr\}\bigg]\text{ with }\tilde{a}_{bi}^*=\frac{\tilde{a}_{bw}^*}{\tilde{a}_{iw}^*},\text{ where } i\in D.\label{5op_PCS_2}
		\end{align}
	\end{subequations}
	From \cite[Theorem 4]{wu2023analytical}, the optimal interval-weight of criterion $c_i$, $i\in C$, is $[{w_i^l}^*,{w_i^u}^*]$, where
	\begin{equation}\label{5opt_weights_1}
		{w_i^l}^*=\frac{\inf\{\tilde{a}_{iw}^*\}}{\inf\{\tilde{a}_{iw}^*\}+\displaystyle\sum_{\substack{j\in C\\j\neq i}}\sup\{\tilde{a}_{jw}^*\}},\quad {w_i^u}^*=\frac{\sup\{\tilde{a}_{iw}^*\}}{\sup\{\tilde{a}_{iw}^*\}+\displaystyle\sum_{\substack{j\in C\\j\neq i}}\inf\{\tilde{a}_{jw}^*\}}.
	\end{equation} 
	Therefore 
	\begin{equation}\label{5opt_weights}
		\begin{split}
			{w_b^l}^*&=\frac{\tilde{a}_{bw}^*}{1+\tilde{a}_{bw}^*+\displaystyle\sum_{j\in D}\min\biggl\{a_{jw}+\epsilon^*,\frac{\tilde{a}_{bw}^*}{a_{bj}-\epsilon^*}\biggr\}},\quad
			{w_b^u}^*=\frac{\tilde{a}_{bw}^*}{1+\tilde{a}_{bw}^*+\displaystyle\sum_{j\in D}\max\biggl\{a_{jw}-\epsilon^*,\frac{\tilde{a}_{bw}^*}{a_{bj}+\epsilon^*}\biggr\}},\\
			{w_w^l}^*&=\frac{1}{1+\tilde{a}_{bw}^*+\displaystyle\sum_{j\in D}\min\biggl\{a_{jw}+\epsilon^*,\frac{\tilde{a}_{bw}^*}{a_{bj}-\epsilon^*}\biggr\}},\quad
			{w_w^u}^*=\frac{1}{1+\tilde{a}_{bw}^*+\displaystyle\sum_{j\in D}\max\biggl\{a_{jw}-\epsilon^*,\frac{\tilde{a}_{bw}^*}{a_{bj}+\epsilon^*}\biggr\}},\\
			{w_i^l}^*&=\frac{\max\biggl\{a_{iw}-\epsilon^*,\frac{\tilde{a}_{bw}^*}{a_{bi}+\epsilon^*}\biggr\}}{1+\tilde{a}_{bw}^*+\max\biggl\{a_{iw}-\epsilon^*,\frac{\tilde{a}_{bw}^*}{a_{bi}+\epsilon^*}\biggr\}+\displaystyle\sum_{\substack{j\in D\\j\neq i}}\min\biggl\{a_{jw}+\epsilon^*,\frac{\tilde{a}_{bw}^*}{a_{bj}-\epsilon^*}\biggr\}},\\
			{w_i^u}^*&=\frac{\min\biggl\{a_{iw}+\epsilon^*,\frac{\tilde{a}_{bw}^*}{a_{bi}-\epsilon^*}\biggr\}}{1+\tilde{a}_{bw}^*+\min\biggl\{a_{iw}+\epsilon^*,\frac{\tilde{a}_{bw}^*}{a_{bi}-\epsilon^*}\biggr\}+\displaystyle\sum_{\substack{j\in D\\j\neq i}}\max\biggl\{a_{jw}-\epsilon^*,\frac{\tilde{a}_{bw}^*}{a_{bj}+\epsilon^*}\biggr\}},\quad \text{where } i\in D.
		\end{split}
	\end{equation}
	\subsubsection{A Secondary Objective Function}
	To identify the optimal weight set that best preserves the given data, Wu et al. \cite{wu2023analytical} proposed the incorporation of a secondary objective function.\\\\
	In this approach, an optimally modified PCS having the minimum value of $\max\{|\tilde{a}_{bi}^*-a_{bi}|,|\tilde{a}_{iw}^*-a_{iw}|\}$ for all $i\in C$ is selected as the best optimally modified PCS. As shown in \cite[Theorem 5]{wu2023analytical}, for $i\in D$, the minimum possible value of $\max\{|\tilde{a}_{bi}^*-a_{bi}|,|\tilde{a}_{iw}^*-a_{iw}|\}$ is the optimal objective value of problem 
	\begin{equation}\label{5optimization_5}
		\begin{split}
			&\min\eta_i \\
			&\text{subject to:}\\	
			&\tilde{a}_{bi}-a_{bi}=\eta_{bi}, \quad \tilde{a}_{iw}-a_{iw}=\eta_{iw},\quad(a_{bi}+\eta_{bi})\times(a_{iw}+\eta_{iw})=\tilde{a}_{bw}^*,\\
			&0\leq \eta_{bi},\eta_{iw}\leq\eta_i
		\end{split}
	\end{equation} 
	if $a_{bi}\times a_{iw}\leq\tilde{a}_{bw}^*$ and the optimal objective value of problem 
	\begin{equation}\label{5optimization_6}
		\begin{split}
			&\min\eta_i \\
			&\text{subject to:}\\	
			&a_{bi}-\tilde{a}_{bi}=\eta_{bi}, \quad a_{iw}-\tilde{a}_{iw}=\eta_{iw},\quad(a_{bi}-\eta_{bi})\times(a_{iw}-\eta_{iw})=\tilde{a}_{bw}^*,\\
			&0\leq \eta_{bi},\eta_{iw}\leq\eta_i
		\end{split}
	\end{equation}
	if $a_{bi}\times a_{iw}>\tilde{a}_{bw}^*$, where $\tilde{a}_{bw}^*$ is as in equation \eqref{5op_PCS_1}. Furthermore, \cite[Theorem 6]{wu2023analytical} states that the PCS 	 
	\begin{equation}\label{5best_pcs}
		\begin{split}
		&\tilde{a}_{bi}^*=\begin{cases}
				a_{bi}+\eta_i^* \quad \text{if } a_{bi}\times a_{iw}\leq \tilde{a}_{bw}^*,\\
				a_{bi}-\eta_i^* \quad \text{if } a_{bi}\times a_{iw}>\tilde{a}_{bw}^*,
			\end{cases}\\
		&\tilde{a}_{iw}^*=\begin{cases}
				a_{iw}+\eta_i^* \quad \text{if } a_{bi}\times a_{iw}\leq\tilde{a}_{bw}^*,\\
				a_{iw}-\eta_i^* \quad \text{if } a_{bi}\times a_{iw}>\tilde{a}_{bw}^*,\text{ where}
			\end{cases}
		\end{split}		
	\end{equation}
	\begin{equation}\label{5CV5}
		\eta_i^*= \bigg|\frac{a_{bi}+a_{iw}-\sqrt{(a_{bi}+a_{iw})^2-4(a_{bi}\times a_{iw}-\tilde{a}_{bw}^*)}}{2}\bigg|,\quad i\in D,
	\end{equation}	
	is the only best optimally modified PCS. Using this PCS, the best optimal weight set is obtained using equation \eqref{5equivalence}.
	\section{Research Gap}
	Despite the wide applicability of the nonlinear BWM, several research gaps remain to be addressed.\\\\
	(i)\quad\textbf{Incompatibility of the existing analytical framework with some scales:}\\\\
		The existing analytical approach is fully compatible with the Saaty scale but proves incompatible with some other scales such as the Salo-H{\"a}m{\"a}l{\"a}inen scale, the Lootsma scale, and the Donegan-Dodd-McMaster scale.\\\\
		\textbf{Example 1:} Let $C=\{c_1,c_2,\ldots,c_5\}$ be the set of decision criteria with $c_1$ as the best and $c_5$ as the worst criterion. Let $A_b=(1,9,3,1.8571,9)$ be the best-to-other vector and $A_w=(9,1.5,4,3,1)^T$ be the other-to-worst vector, obtained by quantifying the preferences given in the form of linguistic terms using the Salo-H{\"a}m{\"a}l{\"a}inen scale.\\\\
		Here, $D_1=\{4\}$, $D_2=\{2,3\}$, and $D_3=\emptyset$. By equation \eqref{5CV}, we get $\epsilon_2=0.4056$, $\epsilon_3=0.3944$, and $\epsilon_4=0.5363$. This gives $i_1=4$ and $i_2=2$. Now, by equation \eqref{5CV'}, we get $\epsilon_{4,2}=0.5163$. Note that $(a_{12}-\epsilon_4)\times (a_{25}-\epsilon_4)<a_{15}-\epsilon_4$. So, equation \eqref{5op_obj} gives $\epsilon^*=\epsilon_4=0.5363$. The optimal interval-weights computed using equation \eqref{5opt_weights} are given in Table \ref{5example_1}.\\\\
		Now, from equation \eqref{5op_PCS_1}, we have $\tilde{a}_{bw}^*=8.4638$. By equation \eqref{5CV5}, we get $\eta_2^*=0.5038$, $\eta_3^*=0.5481$, and $\eta_4^*=0.5363$. From equation \eqref{5best_pcs}, it follows that the best optimally modified PCS is $$\tilde{A}_b^*=(1,8.4962,2.4519,2.3934,8.4638),\quad \tilde{A}_w^*=(8.4638,0.9962,3.4519,3.5363,1)^T.$$ The best optimal weight set calculated using equation \eqref{5equivalence} is presented in Table \ref{5example_1}.\\\\
		It is important to note that the calculated interval-weights are not well-defined, as the lower bound of each interval-weight exceeds the upper bound. Additionally, the value of $\eta_3^*$ is higher than $\epsilon^*$, meaning the resulting best optimally modified PCS should not even be considered as an optimally modified PCS. Consequently, the obtained best optimal weight set should not even be regarded as an optimal weight set.\\\\
		We calculated the actual $\epsilon^*$, optimal interval-weights and the best optimally modified PCS using MATLAB. These $\epsilon^*$ and optimal interval-weights are given in Table \refeq{5example_1}. The best optimally modified PCS is $$\tilde{A}_b^*=(1,8.4999,2.4578,2.3993,8.4987), \quad \tilde{A}_w^*=(8.4987,0.9999,3.4578,3.5422,1)^T.$$ Using equation \eqref{5equivalence}, we derived the best optimal weight set, which is also provided in Table \ref{5example_1}.
		\begin{table}[H]
			\centering
			\caption{Weights and $\epsilon^*$ for Example 1\label{5example_1}}
			\begin{tabular}{@{}ccccccc@{}}
				\toprule[0.1em]
				&\phantom{}&\multicolumn{2}{c}{Analytical approach \cite{wu2023analytical}}&\phantom{}&\multicolumn{2}{c}{Actual value}\\
				\cmidrule{3-4}\cmidrule{6-7}
				Criterion&&Interval-weight&Best weight&&Interval-weight&Best weight\\
				\midrule
				$c_1$&&$[0.4854,0.4857]$&$0.4851$&&$[0.4855,0.4868]$&$0.4857$\\
				$c_2$&&$[0.0554,0.0573]$&$0.0571$&&$[0.0549,0.0574]$&$0.0571$\\
				$c_3$&&$[0.1983,0.1974]$&$0.1978$&&$[0.1975,0.1981]$&$0.1976$\\
				$c_4$&&$[0.2028,0.2029]$&$0.2027$&&$[0.2024,0.2029]$&$0.2024$\\
				$c_5$&&$[0.0573,0.0574]$&$0.0573$&&$[0.0571,0.0573]$&$0.0572$\\
				\midrule
				$\epsilon^*$&&\multicolumn{2}{c}{$0.5363$}&&\multicolumn{2}{c}{$0.5422$}\\
				\bottomrule[0.1em]				
			\end{tabular}
		\end{table}
		\hspace{-0.75cm}
		\textbf{Example 2:} Let $C=\{c_1,c_2,\ldots,c_5\}$ be the set of decision criteria with $c_1$ as the best and $c_5$ as the worst criterion. Let $A_b=(1,16,4\sqrt{2},2\sqrt{2},16)$ be the best-to-other vector and $A_w=(16,2\sqrt{2},4\sqrt{2},\sqrt{2},1)^T$ be the other-to-worst vector, obtained by quantifying the preferences given in the form of linguistic terms using the Lootsma scale.\\\\
		Here, $D_1=\{4\}$, $D_2=\{2,3\}$, and $D_3=\emptyset$. By equation \eqref{5CV}, we get $\epsilon_2=1.6054$, $\epsilon_3=1.4764$, and $\epsilon_4=1.7228$. This gives $i_1=4$ and $i_2=2$. Now, by equation \eqref{5CV'}, we get $\epsilon_{4,2}=1.7882$. Note that $(a_{12}-\epsilon_4)\times (a_{25}-\epsilon_4)>a_{15}-\epsilon_4$ and $(a_{14}+\epsilon_2)\times (a_{45}+\epsilon_2)<a_{15}+\epsilon_2$. So, equation \eqref{5op_obj} gives $\epsilon^*=\epsilon_{4,2}=1.7882$. The optimal interval-weights computed using equation \eqref{5opt_weights} are given in Table \ref{5example_2}.\\\\
		Now, from equation \eqref{5op_PCS_1}, we have $\tilde{a}_{bw}^*=14.7841$. By equation \eqref{5CV5}, we get $\eta_2^*=1.7882$, $\eta_3^*=1.8118$, and $\eta_4^*=1.7882$. From equation \eqref{5best_pcs}, it follows that the best optimally modified PCS is $$\tilde{A}_b^*=(1,14.2118,3.8451,4.6166,14.7841),\quad \tilde{A}_w^*=(14.7841,1.0403,3.8451,3.2024,1)^T.$$ The best optimal weight set calculated using equation \eqref{5equivalence} is presented in Table \ref{5example_2}.\\\\
		It is important to note that the calculated interval-weights are not well-defined, as the lower bound of each interval-weight exceeds the upper bound. Additionally, the value of $\eta_3^*$ is higher than $\epsilon^*$, meaning the resulting best optimally modified PCS should not even be considered an optimally modified PCS. Consequently, the obtained best optimal weight set should not even be regarded as an optimal weight set.\\\\
		We calculated the actual $\epsilon^*$, optimal interval-weights and the best optimally modified PCS using MATLAB. These $\epsilon^*$ and optimal interval-weights are given in Table \refeq{5example_2}. The best optimally modified PCS is $$\tilde{A}_b^*=(1,14.2179,3.8569,4.6283,14.8760), \quad \tilde{A}_w^*=(14.8760,1.0463,3.8569,3.2141,1)^T.$$ Using equation \eqref{5equivalence}, we derived the best optimal weight set, which is also provided in Table \ref{5example_2}.
		\begin{table}[H]
			\centering
			\caption{Weights and $\epsilon^*$ for Example 2\label{5example_2}}
			\begin{tabular}{@{}ccccccc@{}}
				\toprule[0.1em]
				&\phantom{}&\multicolumn{2}{c}{Analytical approach \cite{wu2023analytical}}&\phantom{}&\multicolumn{2}{c}{Actual value}\\
				\cmidrule{3-4}\cmidrule{6-7}
				Criterion&&Interval-weight&Best weight&&Interval-weight&Best weight\\
				\midrule
				$c_1$&&$[0.6199,0.6187]$&$0.6193$&&$[0.6200,0.6205]$&$0.6200$\\
				$c_2$&&$[0.0436,0.0435]$&$0.0436$&&$[0.0429,0.0437]$&$0.0436$\\
				$c_3$&&$[0.1619,0.1602]$&$0.1611$&&$[0.1607,0.1609]$&$0.1607$\\
				$c_4$&&$[0.1343,0.1340]$&$0.1341$&&$[0.1340,0.1341]$&$0.1340$\\
				$c_5$&&$[0.0419,0.0418]$&$0.0419$&&$[0.0416,0.0417]$&$0.0417$\\
				\midrule
				$\epsilon^*$&&\multicolumn{2}{c}{$1.7882$}&&\multicolumn{2}{c}{$1.7999$}\\
				\bottomrule[0.1em]				
			\end{tabular}
		\end{table}
		\hspace{-0.75cm}
		\textbf{Example 3:} Let $C=\{c_1,c_2,\ldots,c_5\}$ be the set of decision criteria with $c_1$ as the best and $c_5$ as the worst criterion. Let $A_b=(1,5.8284,3.2289,1,5.8284)$ be the best-to-other vector and $A_w=(5.8284,1.967,3.2289,1.967,1)^T$ be the other-to-worst vector, obtained by quantifying the preferences given in the form of linguistic terms using the Donegan-Dodd-McMaster scale.\\\\
		Here, $D_1=\{4\}$, $D_2=\{2,3\}$, and $D_3=\emptyset$. By equation \eqref{5CV}, we get $\epsilon_2=0.6959$, $\epsilon_3=0.6781$, and $\epsilon_4=0.8086$. This gives $i_1=4$ and $i_2=2$. Now, by equation \eqref{5CV'}, we get $\epsilon_{4,2}=0.8825$. Note that $(a_{12}-\epsilon_4)\times (a_{25}-\epsilon_4)>a_{15}-\epsilon_4$ and $(a_{14}+\epsilon_2)\times (a_{45}+\epsilon_2)<a_{15}+\epsilon_2$. So, equation \eqref{5op_obj} gives $\epsilon^*=\epsilon_{4,2}=0.8825$. The optimal interval-weights computed using equation \eqref{5opt_weights} are given in Table \ref{5example_3}.\\\\
		Now, from equation \eqref{5op_PCS_1}, we have $\tilde{a}_{bw}^*=5.3640$. By equation \eqref{5CV5}, we get $\eta_2^*=0.8825$, $\eta_3^*=0.9129$, and $\eta_4^*=0.8825$. From equation \eqref{5best_pcs}, it follows that the best optimally modified PCS is $$\tilde{A}_b^*=(1,4.9459,2.316,1.8825,5.3640),\quad \tilde{A}_w^*=(5.3640,1.0845,2.316,2.8495,1)^T.$$ The best optimal weight set calculated using equation \eqref{5equivalence} is presented in Table \ref{5example_3}.\\\\
		It is important to note that the calculated interval-weights are not well-defined, as the lower bound of each interval-weight exceeds the upper bound. Additionally, the value of $\eta_3^*$ is higher than $\epsilon^*$, meaning the resulting best optimally modified PCS should not even be considered an optimally modified PCS. Consequently, the obtained best optimal weight set should not even be regarded as an optimal weight set.\\\\
		We calculated the actual $\epsilon^*$, optimal interval-weights and the best optimally modified PCS using MATLAB. These $\epsilon^*$ and optimal interval-weights are given in Table \refeq{5example_3}. The best optimally modified PCS is $$\tilde{A}_b^*=(1,4.9577,2.3314,1.8975,5.4354), \quad \tilde{A}_w^*=(5.4354,1.0963,2.3314,2.8645,1)^T.$$ Using equation \eqref{5equivalence}, we derived the best optimal weight set, which is also provided in Table \ref{5example_3}.
		\begin{table}[H]
			\centering
			\caption{Weights and $\epsilon^*$ for Example 3\label{5example_3}}
			\begin{tabular}{@{}ccccccc@{}}
				\toprule[0.1em]
				&\phantom{}&\multicolumn{2}{c}{Analytical approach \cite{wu2023analytical}}&\phantom{}&\multicolumn{2}{c}{Actual value}\\
				\cmidrule{3-4}\cmidrule{6-7}
				Criterion&&Interval-weight&Best weight&&Interval-weight&Best weight\\
				\midrule
				$c_1$&&$[0.4263,0.4242]$&$0.4252$&&$[0.4269,0.4280]$&$0.4270$\\
				$c_2$&&$[0.0862,0.0858]$&$0.0860$&&$[0.0842,0.0866]$&$0.0861$\\
				$c_3$&&$[0.1856,0.1817]$&$0.1836$&&$[0.1831,0.1836]$&$0.1832$\\
				$c_4$&&$[0.2264,0.2254]$&$0.2259$&&$[0.2250,0.2255]$&$0.2251$\\
				$c_5$&&$[0.0795,0.0791]$&$0.0793$&&$[0.0785,0.0787]$&$0.0786$\\
				\midrule
				$\epsilon^*$&&\multicolumn{2}{c}{$0.8825$}&&\multicolumn{2}{c}{$0.8975$}\\
				\bottomrule[0.1em]				
			\end{tabular}
		\end{table}
		\hspace{-0.7cm}
		(ii)\quad\textbf{Incompatibility of the existing framework in the presence of multiple DMs:}\\\\
		The framework does not work even with the Saaty scale in group decision contexts where the Aggregation of Individual Judgements (AIJ) approach is employed \cite{brunelli2014introduction}.\\\\
		\textbf{Example 4:} Let $C=\{c_1,c_2,\ldots,c_5\}$ be the set of decision criteria with $c_1$ as the best and $c_5$ as the worst criterion. Let $E_1$ and $E_2$ be two homogeneous DMs. Let $(A_b)_1=(1,5,1,3,7)$ be the best-to-other vector and $(A_w)_1=(7,2,7,1,1)^T$ be the other-to-worst vector for $E_1$, and let $(A_b)_2=(1,2,5,2,7)$ be the best-to-other vector and $(A_w)_2=(7,5,3,3,1)^T$ be the other-to-worst vector for $E_2$. These vectors are obtained by quantifying the preferences given in the form of linguistic terms using the Saaty scale.\\\\
		Now, we aggregate these individual judgments using the geometric mean method. Thus, we get $A_b=(1,\sqrt{10},\sqrt{5},\sqrt{6},7)$ as the aggregated best-to-other vector and $A_w=(7,\sqrt{10},\sqrt{21},\sqrt{3},1)^T$ as the aggregated other-to-worst vector.\\\\
		Here, $D_1=\{4\}$, $D_2=\{2,3\}$, and $D_3=\emptyset$. By equation \eqref{5CV}, we get $\epsilon_2=0.4355$, $\epsilon_3=0.4401$, and $\epsilon_4=0.4865$. This gives $i_1=4$ and $i_2=3$. Now, by equation \eqref{5CV'}, we get $\epsilon_{4,3}=0.5458$. Note that $(a_{13}-\epsilon_4)\times (a_{35}-\epsilon_4)>a_{15}-\epsilon_4$ and $(a_{14}+\epsilon_3)\times (a_{45}+\epsilon_3)<a_{15}+\epsilon_3$. So, equation \eqref{5op_obj} gives $\epsilon^*=\epsilon_{4,3}=0.5458$. The optimal interval-weights computed using equation \eqref{5opt_weights} are given in Table \ref{5example_4}.\\\\
		Now, from equation \eqref{5op_PCS_1}, we have $\tilde{a}_{bw}^*=6.8230$. By equation \eqref{5CV5}, we get $\eta_2^*=0.5502$, $\eta_3^*=0.5458$, and $\eta_4^*=0.5458$. From equation \eqref{5best_pcs}, it follows that the best optimally modified PCS is $$\tilde{A}_b^*=(1,2.6121,1.6902,2.9953,6.8230),\quad \tilde{A}_w^*=(6.8230,2.6121,4.0367,2.2779,1)^T.$$ The best optimal weight set calculated using equation \eqref{5equivalence} is presented in Table \ref{5example_4}.\\\\
		It is important to note that the calculated interval-weights are not well-defined, as the lower bound of each interval-weight exceeds the upper bound. Additionally, the value of $\eta_2^*$ is higher than $\epsilon^*$, meaning the resulting best optimally modified PCS should not even be considered an optimally modified PCS. Consequently, the obtained best optimal weight set should not even be regarded as an optimal weight set.\\\\
		We calculated the actual $\epsilon^*$, optimal interval-weights and the best optimally modified PCS using MATLAB. These $\epsilon^*$ and optimal interval-weights are given in Table \refeq{5example_4}. The best optimally modified PCS is $$\tilde{A}_b^*=(1,2.6143,1.6923,2.9975,6.8344), \quad \tilde{A}_w^*=(6.8344,2.6143,4.0388,2.2801,1)^T.$$ Using equation \eqref{5equivalence}, we derived the best optimal weight set, which is also provided in Table \ref{5example_4}.
		\begin{table}[H]
			\centering
			\caption{Weights and $\epsilon^*$ for Example 4\label{5example_4}}
			\begin{tabular}{@{}ccccccc@{}}
				\toprule[0.1em]
				&\phantom{}&\multicolumn{2}{c}{Analytical approach \cite{wu2023analytical}}&\phantom{}&\multicolumn{2}{c}{Actual value}\\
				\cmidrule{3-4}\cmidrule{6-7}
				Criterion&&Interval-weight&Best weight&&Interval-weight&Best weight\\
				\midrule
				$c_1$&&$[0.4099,0.4047]$&$0.4074$&&$[0.4074,0.4077]$&$0.4076$\\
				$c_2$&&$[0.1615,0.1506]$&$0.1559$&&$[0.1558,0.1560]$&$0.1559$\\
				$c_3$&&$[0.2425,0.2394]$&$0.2410$&&$[0.2407,0.2413]$&$0.2409$\\
				$c_4$&&$[0.1369,0.1351]$&$0.1360$&&$[0.1359,0.1360]$&$0.1360$\\
				$c_5$&&$[0.0601,0.0593]$&$0.0597$&&$[0.0596,0.0597]$&$0.0596$\\
				\midrule
				$\epsilon^*$&&\multicolumn{2}{c}{$0.5458$}&&\multicolumn{2}{c}{$0.5480$}\\
				\bottomrule[0.1em]				
			\end{tabular}
		\end{table}
		\hspace{-0.7cm}
		(iii)\quad\textbf{Non-well-defined values of CI for the Saaty scale:}\\\\
		The values of CI for the Saaty scale, presented in Table \ref{5ci}, do not serve as an upper bound for the set of $\epsilon^*$ corresponding to the PCSs having the given $a_{bw}$, and thus, is not well-defined.\\\\
		\textbf{Example 5:} Let $\{c_1,c_2,c_3,c_4\}$ be the set of decision criteria with $c_1$ as the best and $c_4$ as the worst criterion. Let $A_b=(1,1,2,2)$ be the best-to-other vector and $A_w=(2,1,2,1)^T$ be the other-to-worst vector, obtained by quantifying the preferences given in the form of linguistic terms using the Saaty scale.\\\\
		Here, we get a unique optimal weight set, which along with $\epsilon^*$, is given in Table \ref{5example_5}.\\\\
		Note that $\epsilon^*=0.5>0.4384= \text{CI}_{a_{bw}=2}$. This implies that the values of CI for the Saaty scale given in Table \ref{5ci} are not well-defined.
		\begin{table}[H]
			\centering
			\caption{Weights and $\epsilon^*$ for Example 5\label{5example_5}}
			\begin{tabular}{ccc}
				\toprule[0.1em]
				Criterion&\phantom{}&Weight\\
				\midrule
				$c_1$&&$0.36$\\
				$c_2$&&$0.24$\\
				$c_3$&&$0.24$\\
				$c_4$&&$0.16$\\
				\midrule
				$\epsilon^*$&&$0.5$\\
				\bottomrule[0.1em]				
			\end{tabular}
		\end{table}
		\hspace{-0.75cm}
		(iv)\quad\textbf{Weight differences among criteria with equal preference in the presence of multiple best/worst criteria:}\\\\
		In instances where multiple criteria are equally qualified as the best (or worst), the conventional approach involves arbitrarily selecting one as the best (or worst) for computational purposes. However, the chosen one might get a different weight than the others even though they are of equal preference.\\\\
		\textbf{Example 6:} Let $\{c_1,c_2,c_3,c_4\}$ be the set of decision criteria with $c_1$ and $c_2$ as the best criteria and $c_4$ as the worst criterion. Here, we select $c_1$ as the best criterion and proceed accordingly. Let $A_b=(1,1,2,7)$ be the best-to-other vector and $A_w=(7,7,3,1)^T$ be the other-to-worst vector, obtained by quantifying the preferences given in the form of linguistic terms using the Saaty scale.\\\\
		The interval-weights, the best optimal weight set and $\epsilon^*$ are given in Table \ref{5example_6}.\\\\
		Note that for $c_1$ and $c_2$, neither the interval-weights nor the best weights coincide, despite both criteria being of equal preferences.
		\begin{table}[H]
			\centering
			\caption{Weights and $\epsilon^*$ for Example 6\label{5example_6}}
			\begin{tabular}{@{}cccc@{}}
				\toprule[0.1em]
				Criterion&\phantom{}&Interval-weight&Best weight\\
				\midrule
				$c_1$&&$[0.3765,0.3833]$&$0.3803$\\
				$c_2$&&$[0.3833,0.3944]$&$0.3882$\\
				$c_3$&&$[0.1741,0.1773]$&$0.1759$\\
				$c_4$&&$[0.0551,0.0561]$&$0.0556$\\
				\midrule
				$\epsilon^*$&&\multicolumn{2}{c}{$0.1623$}\\
				\bottomrule[0.1em]				
			\end{tabular}
		\end{table}
\section{A Generalized Analytical Framework for the Nonlinear BWM}
In this section, we propose an analytical approach compatible with any scale and any number of DMs. We also modify the original optimization model to ensure equal weights for all best and worst criteria in cases where multiple best or worst criteria are present. Furthermore, we derive formulas for CI and CR.
\subsection{Calculation of Weights}
This subsection consists of two parts: weight calculation when a unique best and worst criterion exists, and weight calculation for cases with multiple best or worst criteria.
\subsubsection{Unique Best and Worst Criterion}
It is easy to verify that equations \eqref{5op_PCS_1} and \eqref{5opt_weights} hold in general, i.e., they are valid for any scale and any number of DMs. Therefore, to determine the optimal interval-weights, it suffices to derive the general expressions for $\epsilon^*$ and $\tilde{a}_{bw}^*$. In the following, our objective is to compute these values. We also derive the best optimally modified PCS that leads to the best optimal weight set.
\begin{proposition}\label{5lower}
	Let $\epsilon_i$ and $\epsilon_{i,j}$ be as in equations \eqref{5CV} and \eqref{5CV'} respectively, and let $\eta^*$ be the optimal objective value of problem \eqref{5optimization_2}. Then $\epsilon_i\leq\eta^*$ and $\epsilon_{i,j}\leq\eta^*$ for all $i,j\in D$.
\end{proposition}
\begin{proof}
	Let $(\tilde{A}_b^*,\tilde{A}_w^*)$ be an optimally modified PCS. Let $|\tilde{a}_{bi}^*-a_{bi}|=\eta_{bi}$, $|\tilde{a}_{iw}^*-a_{iw}|=\eta_{iw}$ and $|\tilde{a}_{bw}^*-a_{bw}|=\eta_{bw}$. Therefore, $0\leq\eta_{bi},\eta_{iw},\eta_{bw}\leq \eta^*$. Also, $\tilde{a}_{bi}^*\in\{a_{bi}+\eta_{bi}, a_{bi}-\eta_{bi}\}$, $\tilde{a}_{iw}^*\in\{a_{iw}+\eta_{iw}, a_{iw}-\eta_{iw}\}$ and $\tilde{a}_{bw}^*\in\{a_{bw}+\eta_{bw}, a_{bw}-\eta_{bw}\}$.\\\\
	Fix $i\in D$. If $i\in D_3$, then by equation \eqref{5CV}, $\epsilon_i=0$ and we are done. If $i\in D_1$, then by equation \eqref{5CV1}, we have $(a_{bi}+\epsilon_i)\times(a_{iw}+\epsilon_i)=a_{bw}-\epsilon_i$. Now, to prove $\epsilon_i\leq \eta^*$, it suffices to show that at least one of the inequalities $\epsilon_i\leq\eta_{bi}$, $\epsilon_i\leq\eta_{iw}$ or $\epsilon_i\leq\eta_{bw}$ holds. Suppose, if possible, neither of these inequalities hold. Then we get $a_{bi}+\epsilon_i>\tilde{a}_{bi}^*$, $a_{iw}+\epsilon_i>\tilde{a}_{iw}^*$ and $a_{bw}-\epsilon_i<\tilde{a}_{bw}^*$. This implies $\tilde{a}_{bi}^*\times \tilde{a}_{iw}^*<\tilde{a}_{bw}^*$, which is contradiction as $(\tilde{A}_b^*,\tilde{A}_w^*)$ is consistent. Thus, $\epsilon_i\leq \eta^*$. If $i\in D_2$, then the result follows by applying the same argument as above.\\\\
	Fix $i,j\in D$. If $a_{bi}\times a_{iw}=a_{bj}\times a_{jw}$. then by equation \eqref{5CV'}, $\epsilon_{i,j}=0$ and we are done. If $a_{bi}\times a_{iw}<a_{bj}\times a_{jw}$, then by equation \eqref{5CV3}, we have $(a_{bi}+\epsilon_{i,j})\times (a_{iw}+\epsilon_{i,j})= (a_{bj}-\epsilon_{i,j})\times(a_{jw}-\epsilon_{i,j})$. Now, to prove $\epsilon_{i,j}\leq \eta^*$, it suffices to show that at least one of the inequalities $\epsilon_{i,j}\leq\eta_{bi}$, $\epsilon_{i,j}\leq\eta_{iw}$, $\epsilon_{i,j}\leq\eta_{bj}$ or $\epsilon_{i,j}\leq\eta_{jw}$ holds. Suppose, if possible, neither of these inequalities hold. Then we get $a_{bi}+\epsilon_{i,j}>\tilde{a}_{bi}^*$, $a_{iw}+\epsilon_{i,j}>\tilde{a}_{iw}^*$, $a_{bj}-\epsilon_{i,j}<\tilde{a}_{bj}^*$ and $a_{jw}-\epsilon_{i,j}<\tilde{a}_{jw}^*$. This implies $\tilde{a}_{bi}^*\times \tilde{a}_{iw}^*<\tilde{a}_{bj}^*\times \tilde{a}_{jw}^*$, which is contradiction as $(\tilde{A}_b^*,\tilde{A}_w^*)$ is consistent. Thus, $\epsilon_{i,j}\leq \eta^*$. Hence the proof.
\end{proof}
\hspace{-0.75cm}
Proposition \ref{5lower} establishes that $\max\{\epsilon_i,\epsilon_{i,j}:i,j\in D\}$ serves as a lower bound for $\eta^*$ (and consequently, for $\epsilon^*$).
\begin{theorem}\label{5exact_obj}
	Let $\epsilon_i$ and $\epsilon_{i,j}$ be as in equations \eqref{5CV} and \eqref{5CV'} respectively, and let $\eta^*$ be the optimal objective value of problem \eqref{5optimization_2}. Then the following statements hold.
	\begin{enumerate}
		\item If $\epsilon_{i_0}=\max\{\epsilon_i,\epsilon_{i,j}:i,j\in D\}$ for some $i_0\in D_1$, then $\eta^*=\epsilon_{i_0}$. Also, $(\tilde{A}_b,\tilde{A}_w)$ defined as 
		\begin{equation}\label{5best_pcs_1}
			\begin{split}
			&\tilde{a}_{bi}=\frac{a_{bi}-a_{iw}+\sqrt{(a_{bi}+a_{iw})^2-4\times a_{bi}\times a_{iw}+4(a_{bw}-\epsilon_{i_0})}}{2},\\
			&\tilde{a}_{iw}=\frac{-a_{bi}+a_{iw}+\sqrt{(a_{bi}+a_{iw})^2-4\times a_{bi}\times a_{iw}+4(a_{bw}-\epsilon_{i_0})}}{2}\quad\text{for all } i\in D,\\
			&\tilde{a}_{bw}=a_{bw}-\epsilon_{i_0}
			\end{split}
		\end{equation}
		is an optimally modified PCS.
		\item If $\epsilon_{j_0}=\max\{\epsilon_i,\epsilon_{i,j}:i,j\in D\}$ for some $j_0\in D_2$, then $\eta^*=\epsilon_{j_0}$. Also, $(\tilde{A}_b,\tilde{A}_w)$ defined as 
		\begin{equation}\label{5best_pcs_2}
			\begin{split}
				&\tilde{a}_{bi}=\frac{a_{bi}-a_{iw}+\sqrt{(a_{bi}+a_{iw})^2-4\times a_{bi}\times a_{iw}+4(a_{bw}+\epsilon_{j_0})}}{2},\\
				&\tilde{a}_{iw}=\frac{-a_{bi}+a_{iw}+\sqrt{(a_{bi}+a_{iw})^2-4\times a_{bi}\times a_{iw}+4(a_{bw}+\epsilon_{j_0})}}{2}\quad\text{for all } i\in D,\\
				&\tilde{a}_{bw}=a_{bw}+\epsilon_{j_0}
			\end{split}
		\end{equation}
		is an optimally modified PCS.
		\item If $\epsilon_{i_0,j_0}=\max\{\epsilon_i,\epsilon_{i,j}:i,j\in D\}$ for some $i_0,j_0\in D$, then $\eta^*=\epsilon_{i_0,j_0}$. Also, $(\tilde{A}_b,\tilde{A}_w)$ defined as 
		\begin{equation}\label{5best_pcs_3}
			\begin{split}
				&\tilde{a}_{bi}=\frac{a_{bi}-a_{iw}+\sqrt{(a_{bi}+a_{iw})^2-4\times a_{bi}\times a_{iw}+4(a_{bi_0}+\epsilon_{i_0,j_0})\times (a_{i_0w}+\epsilon_{i_0,j_0})}}{2},\\
				&\tilde{a}_{iw}=\frac{-a_{bi}+a_{iw}+\sqrt{(a_{bi}+a_{iw})^2-4\times a_{bi}\times a_{iw}+4(a_{bi_0}+\epsilon_{i_0,j_0})\times (a_{i_0w}+\epsilon_{i_0,j_0})}}{2}\\
				&\quad\quad\quad\quad\quad\quad\quad\quad\quad\quad\quad\quad\quad\quad\quad\quad\quad\quad\quad\quad\quad\quad\quad\quad\quad\quad\quad\quad\quad\text{for all } i\in D,\\
				&\tilde{a}_{bw}=(a_{bi_0}+\epsilon_{i_0,j_0})\times (a_{i_0w}+\epsilon_{i_0,j_0})
			\end{split}
		\end{equation}
	is an optimally modified PCS.
	\end{enumerate}
	\begin{proof}
		First, assume that $\epsilon_{i_0}=\max\{\epsilon_i,\epsilon_{i,j}:i,j\in D\}$ for some $i_0\in D_1$. It is easy to verify that $(\tilde{A}_b,\tilde{A}_w)$ given in equation \eqref{5best_pcs_1} is consistent. Now, if we prove that
		\begin{equation*}
			|\tilde{a}_{bi}-a_{bi}|\leq \epsilon_{i_0},\quad |\tilde{a}_{iw}-a_{iw}|\leq \epsilon_{i_0},\quad |\tilde{a}_{bw}-a_{bw}|\leq \epsilon_{i_0}
		\end{equation*}
		for all $i\in D$, then it will imply that $\eta^*\leq \epsilon_{i_0}$. This, along with Proposition \ref{5lower}, gives $\eta^*=\epsilon_{i_0}$, and consequently, $(\tilde{A}_b,\tilde{A}_w)$ is an optimally modified PCS.\\\\
		Now, $|\tilde{a}_{bw}-a_{bw}|=|-\epsilon_{i_0}|=\epsilon_{i_0}$. Also, $$|\tilde{a}_{bi}-a_{bi}|=|\tilde{a}_{iw}-a_{iw}|=\bigg|\frac{-(a_{bi}+a_{iw})+\sqrt{(a_{bi}+a_{iw})^2-4\times a_{bi}\times a_{iw}+4(a_{bw}-\epsilon_{i_0})}}{2}\bigg|.$$ Fix $i\in D$. Then there are two possibilities.
		\begin{enumerate}[(i)]
			\item $a_{bi}\times a_{iw}\leq a_{bw}-\epsilon_{i_0}$\\\\
				Here, we get $(a_{bi}+a_{iw})\leq \sqrt{(a_{bi}+a_{iw})^2-4\times a_{bi}\times a_{iw}+4(a_{bw}-\epsilon_{i_0})}$. Therefore, $$\frac{-(a_{bi}+a_{iw})+\sqrt{(a_{bi}+a_{iw})^2-4\times a_{bi}\times a_{iw}+4(a_{bw}-\epsilon_{i_0})}}{2}\geq 0,$$ and thus, 
				\begin{equation}\label{5opt_dev_1}
					|\tilde{a}_{bi}-a_{bi}|=|\tilde{a}_{iw}-a_{iw}|=\frac{-(a_{bi}+a_{iw})+\sqrt{(a_{bi}+a_{iw})^2-4\times a_{bi}\times a_{iw}+4(a_{bw}-\epsilon_{i_0})}}{2}.
				\end{equation} 
				Note that
				\begin{equation}\label{5dev_eq_1}
					\begin{split}
					\bigg(a_{bi}&+\frac{-(a_{bi}+a_{iw})+\sqrt{(a_{bi}+a_{iw})^2-4\times a_{bi}\times a_{iw}+4(a_{bw}-\epsilon_{i_0})}}{2}\bigg)\\
					&\times \bigg(a_{iw}+\frac{-(a_{bi}+a_{iw})+\sqrt{(a_{bi}+a_{iw})^2-4\times a_{bi}\times a_{iw}+4(a_{bw}-\epsilon_{i_0})}}{2}\bigg)=a_{bw}-\epsilon_{i_0}.
					\end{split}
				\end{equation}
				Since $a_{bi}\times a_{iw}\leq a_{bw}-\epsilon_{i_0}< a_{bw}$, by equation \eqref{5CV1}, we have $(a_{bi}+\epsilon_i)\times(a_{iw}+\epsilon_i)=a_{bw}-\epsilon_i$. Now, $\epsilon_i\leq \epsilon_{i_0}$ implies that $a_{bw}-\epsilon_{i_0}\leq a_{bw}-\epsilon_i$, i.e.,
				\begin{equation*}
					\begin{split}
						\bigg(a_{bi}&+\frac{-(a_{bi}+a_{iw})+\sqrt{(a_{bi}+a_{iw})^2-4\times a_{bi}\times a_{iw}+4(a_{bw}-\epsilon_{i_0})}}{2}\bigg)\\
						&\times \bigg(a_{iw}+\frac{-(a_{bi}+a_{iw})+\sqrt{(a_{bi}+a_{iw})^2-4\times a_{bi}\times a_{iw}+4(a_{bw}-\epsilon_{i_0})}}{2}\bigg)\\
						&\quad\quad\quad\quad\quad\quad\quad\quad\quad\quad\quad\quad\quad\quad\quad\quad\quad\quad\quad\quad\quad\quad\quad\quad\quad\quad\leq (a_{bi}+\epsilon_i)\times(a_{iw}+\epsilon_i).
					\end{split}
				\end{equation*}
				This gives $$|\tilde{a}_{bi}-a_{bi}|=\frac{-(a_{bi}+a_{iw})+\sqrt{(a_{bi}+a_{iw})^2-4\times a_{bi}\times a_{iw}+4(a_{bw}-\epsilon_{i_0})}}{2}\leq \epsilon_i\leq \epsilon_{i_0}.$$
				Similarly, we get $|\tilde{a}_{bi}-a_{bi}|\leq \epsilon_{i_0}$.\\
			\item $a_{bi}\times a_{iw}> a_{bw}-\epsilon_{i_0}$\\\\
			Here, we get $\sqrt{(a_{bi}+a_{iw})^2-4\times a_{bi}\times a_{iw}+4(a_{bw}-\epsilon_{i_0})}<a_{bi}+a_{iw}$. Therefore, $$\frac{a_{bi}+a_{iw}-\sqrt{(a_{bi}+a_{iw})^2-4\times a_{bi}\times a_{iw}+4(a_{bw}-\epsilon_{i_0})}}{2}\geq 0,$$ and thus,
			\begin{equation}\label{5opt_dev_2}
				|\tilde{a}_{bi}-a_{bi}|=|\tilde{a}_{iw}-a_{iw}|=\frac{a_{bi}+a_{iw}-\sqrt{(a_{bi}+a_{iw})^2-4\times a_{bi}\times a_{iw}+4(a_{bw}-\epsilon_{i_0})}}{2}.
			\end{equation}
			Note that
			\begin{equation*}
				\begin{split}
					\bigg(a_{bi}&-\frac{a_{bi}+a_{iw}-\sqrt{(a_{bi}+a_{iw})^2-4\times a_{bi}\times a_{iw}+4(a_{bw}-\epsilon_{i_0})}}{2}\bigg)\\
					&\times \bigg(a_{iw}-\frac{a_{bi}+a_{iw}-\sqrt{(a_{bi}+a_{iw})^2-4\times a_{bi}\times a_{iw}+4(a_{bw}-\epsilon_{i_0})}}{2}\bigg)=a_{bw}-\epsilon_{i_0}.
				\end{split}
			\end{equation*}
			By equation \eqref{5CV1}, we have $(a_{bi_0}+\epsilon_{i_0})\times(a_{i_0w}+\epsilon_{i_0})=a_{bw}-\epsilon_{i_0}$. Since $a_{bi}\times a_{iw}>a_{bw}-\epsilon_{i_0}$, we get $a_{bi_0}\times a_{i_0w}<a_{bi}\times a_{iw}$. So, by equation \eqref{5CV3}, we have $(a_{bi_0}+\epsilon_{i,i_0})\times(a_{i_0w}+\epsilon_{i,i_0})=(a_{bi}-\epsilon_{i,i_0})\times (a_{iw}-\epsilon_{i,i_0})$. Now, $\epsilon_{i,i_0}\leq\epsilon_{i_0}$ implies that $(a_{bi}-\epsilon_{i,i_0})\times (a_{iw}-\epsilon_{i,i_0})\leq a_{bw}-\epsilon_{i_0}$, i.e.,
			\begin{equation*}
				\begin{split}
					&(a_{bi}-\epsilon_{i,i_0})\times (a_{iw}-\epsilon_{i,i_0})\\
					&\quad\quad\quad\quad\quad\leq \bigg(a_{bi}-\frac{a_{bi}+a_{iw}-\sqrt{(a_{bi}+a_{iw})^2-4\times a_{bi}\times a_{iw}+4(a_{bw}-\epsilon_{i_0})}}{2}\bigg)\\
					&\quad\quad\quad\quad\quad\quad\quad\quad\times \bigg(a_{iw}-\frac{a_{bi}+a_{iw}-\sqrt{(a_{bi}+a_{iw})^2-4\times a_{bi}\times a_{iw}+4(a_{bw}-\epsilon_{i_0})}}{2}\bigg).
				\end{split}
			\end{equation*}
			This gives $$|\tilde{a}_{bi}-a_{bi}|=\frac{a_{bi}+a_{iw}-\sqrt{(a_{bi}+a_{iw})^2-4\times a_{bi}\times a_{iw}+4(a_{bw}-\epsilon_{i_0})}}{2}\leq \epsilon_{i,i_0}\leq \epsilon_{i_0}.$$
			Similarly, we get $|\tilde{a}_{bi}-a_{bi}|\leq \epsilon_{i_0}$. 
		\end{enumerate}
		The result can be proven using a similar argument as above if $\epsilon_{j_0}=\max\{\epsilon_i,\epsilon_{i,j}:i,j\in D\}$ for some $j_0\in D_2$.\\\\
		Now, assume that $\epsilon_{i_0,j_0}=\max\{\epsilon_i,\epsilon_{i,j}:i,j\in D\}$ for some $i,j\in D$. It is easy to check that $(\tilde{A}_b,\tilde{A}_w)$ given in equation \eqref{5best_pcs_1} is consistent. Now, it is sufficient to prove that
		\begin{equation*}
			|\tilde{a}_{bi}-a_{bi}|\leq \epsilon_{i_0,j_0},\quad |\tilde{a}_{iw}-a_{iw}|\leq \epsilon_{i_0,j_0},\quad |\tilde{a}_{bw}-a_{bw}|\leq \epsilon_{i_0,j_0}
		\end{equation*}
		for all $i\in D$.\\\\
		Let $\tilde{a}_{bw}-a_{bw}=\zeta$. Then there are two possibilities.
		\begin{enumerate}[(i)]
			\item $\zeta\leq 0$\\
			Here, we get $\tilde{a}_{bw}\leq a_{bw}$. Therefore, $a_{bi_0}\times a_{i_0w}<a_{bw}$. So, by equation \eqref{5CV1}, we have $(a_{bi_0}+\epsilon_{i_0})\times(a_{i_0w}+\epsilon_{i_0})=a_{bw}-\epsilon_{i_0}$. Also, $(a_{bi_0}+\epsilon_{i_0,j_0})\times(a_{i_0w}+\epsilon_{i_0,j_0})=\tilde{a}_{bw}=a_{bw}+\zeta$. Now, $\epsilon_{i_0}\leq \epsilon_{i_0,j_0}$ implies that $(a_{bi_0}+\epsilon_{i_0})\times(a_{i_0w}+\epsilon_{i_0}) \leq (a_{bi_0}+\epsilon_{i_0,j_0})\times(a_{i_0w}+\epsilon_{i_0,j_0})$, i.e., $a_{bw}-\epsilon_{i_0}\leq a_{bw}+\zeta$. This gives $|\tilde{a}_{bw}-a_{bw}|=-\zeta\leq \epsilon_{i_0}\leq \epsilon_{i_0,j_0}$.
			\item $\zeta>0$\\
			Here, we get $a_{bw}<\tilde{a}_{bw}$. From equation \eqref{5CV3}, it follows that $\tilde{a}_{bw}=a_{bw}+\zeta=(a_{bj_0}-\epsilon_{i_0,j_0})\times(a_{j_0w}-\epsilon_{i_0,j_0})$. Therefore, $a_{bj_0}\times a_{j_0w}>a_{bw}$. So, by equation \eqref{5CV2}, we have $(a_{bj_0}-\epsilon_{j_0})\times(a_{j_0w}-\epsilon_{j_0})=a_{bw}+\epsilon_{j_0}$. Now, $\epsilon_{j_0}\leq \epsilon_{i_0,j_0}$ implies that $(a_{bj_0}-\epsilon_{i_0,j_0})\times(a_{j_0w}-\epsilon_{i_0,j_0}) \leq (a_{bj_0}-\epsilon_{j_0})\times(a_{j_0w}-\epsilon_{j_0})$, i.e., $a_{bw}+\zeta \leq a_{bw}+\epsilon_{j_0}$. This gives $|\tilde{a}_{bw}-a_{bw}|=\zeta\leq \epsilon_{j_0}\leq \epsilon_{i_0,j_0}$.
		\end{enumerate}
		Now, 
		\begin{equation*}
			\begin{split}
				|\tilde{a}_{bi}-a_{bi}|&=|\tilde{a}_{iw}-a_{iw}|\\
				&=\bigg|\frac{-(a_{bi}+a_{iw})+\sqrt{(a_{bi}+a_{iw})^2-4\times a_{bi}\times a_{iw}+4(a_{bi_0}+\epsilon_{i_0,j_0})\times (a_{i_0w}+\epsilon_{i_0,j_0})}}{2}\bigg|.
			\end{split}
		\end{equation*}
		Fix $i\in D$. Then there are two possibilities.
		\begin{enumerate}[(i)]
			\item $a_{bi}\times a_{iw}\leq (a_{bi_0}+\epsilon_{i_0,j_0})\times(a_{i_0w}+\epsilon_{i_0,j_0})$\\\\
			Here, we get $(a_{bi}+a_{iw})\leq \sqrt{(a_{bi}+a_{iw})^2-4\times a_{bi}\times a_{iw}+4(a_{bi_0}+\epsilon_{i_0,j_0})\times (a_{i_0w}+\epsilon_{i_0,j_0})}$. So, $$\frac{-(a_{bi}+a_{iw})+\sqrt{(a_{bi}+a_{iw})^2-4\times a_{bi}\times a_{iw}+4(a_{bi_0}+\epsilon_{i_0,j_0})\times (a_{i_0w}+\epsilon_{i_0,j_0})}}{2}\geq 0,$$ and thus, 
			\begin{equation}\label{5opt_dev_3}
				\begin{split}
					|\tilde{a}_{bi}-a_{bi}|&=|\tilde{a}_{iw}-a_{iw}|\\
					&=\frac{-(a_{bi}+a_{iw})+\sqrt{(a_{bi}+a_{iw})^2-4\times a_{bi}\times a_{iw}+4(a_{bi_0}+\epsilon_{i_0,j_0})\times (a_{i_0w}+\epsilon_{i_0,j_0})}}{2}.
				\end{split}
			\end{equation} Note that
			\begin{equation}\label{5dev_eq_3}
				\begin{split}
					&\bigg(a_{bi}+\frac{-(a_{bi}+a_{iw})+\sqrt{(a_{bi}+a_{iw})^2-4\times a_{bi}\times a_{iw}+4(a_{bi_0}+\epsilon_{i_0,j_0})\times (a_{i_0w}+\epsilon_{i_0,j_0})}}{2}\bigg)\\
					&\times \bigg(a_{iw}+\frac{-(a_{bi}+a_{iw})+\sqrt{(a_{bi}+a_{iw})^2-4\times a_{bi}\times a_{iw}+4(a_{bi_0}+\epsilon_{i_0,j_0})\times (a_{i_0w}+\epsilon_{i_0,j_0})}}{2}\bigg)\\
					&\quad\quad\quad\quad\quad\quad\quad\quad\quad\quad\quad\quad\quad\quad\quad\quad\quad\quad\quad\quad\quad\quad\quad\quad\quad=(a_{bi_0}+\epsilon_{i_0,j_0})\times(a_{i_0w}+\epsilon_{i_0,j_0}).
				\end{split}
			\end{equation}
			Since $a_{bi}\times a_{iw}\leq (a_{bi_0}+\epsilon_{i_0,j_0})\times(a_{i_0w}+\epsilon_{i_0,j_0})$ and by equation \eqref{5CV3}, we have $(a_{bi_0}+\epsilon_{i_0,j_0})\times(a_{i_0w}+\epsilon_{i_0,j_0})=(a_{bj_0}-\epsilon_{i_0,j_0})\times(a_{j_0w}-\epsilon_{i_0,j_0})$, we get $a_{bi}\times a_{iw}<a_{bj_0}\times a_{j_0w}$. So, from equation \eqref{5CV3}, we have $(a_{bi}+\epsilon_{i,j_0})\times(a_{iw}+\epsilon_{i,j_0})=(a_{bj_0}-\epsilon_{i,j_0})\times(a_{j_0w}-\epsilon_{i,j_0})$. Now, $\epsilon_{i,j_0}\leq\epsilon_{i_0,j_0}$ implies that $(a_{bi_0}+\epsilon_{i_0,j_0})\times(a_{i_0w}+\epsilon_{i_0,j_0}) \leq (a_{bi}+\epsilon_{i,j_0})\times(a_{iw}+\epsilon_{i,j_0})$, i.e.,
			\begin{equation*}
				\begin{split}
					&\bigg(a_{bi}+\frac{-(a_{bi}+a_{iw})+\sqrt{(a_{bi}+a_{iw})^2-4\times a_{bi}\times a_{iw}+4(a_{bi_0}+\epsilon_{i_0,j_0})\times (a_{i_0w}+\epsilon_{i_0,j_0})}}{2}\bigg)\\
					&\times \bigg(a_{iw}+\frac{-(a_{bi}+a_{iw})+\sqrt{(a_{bi}+a_{iw})^2-4\times a_{bi}\times a_{iw}+4(a_{bi_0}+\epsilon_{i_0,j_0})\times (a_{i_0w}+\epsilon_{i_0,j_0})}}{2}\bigg)\\
					&\quad\quad\quad\quad\quad\quad\quad\quad\quad\quad\quad\quad\quad\quad\quad\quad\quad\quad\quad\quad\quad\quad\quad\quad\quad\leq(a_{bi}+\epsilon_{i,j_0})\times(a_{iw}+\epsilon_{i,j_0}).
				\end{split}
			\end{equation*}
			This gives 
			\begin{equation*}
				\begin{split}
					|\tilde{a}_{bi}-a_{bi}|&=\frac{-(a_{bi}+a_{iw})+\sqrt{(a_{bi}+a_{iw})^2-4\times a_{bi}\times a_{iw}+4(a_{bi_0}+\epsilon_{i_0,j_0})\times (a_{i_0w}+\epsilon_{i_0,j_0})}}{2}\\
					&\leq \epsilon_{i,j_0}\leq\epsilon_{i_0,j_0}.
				\end{split}
			\end{equation*}
			Similarly, we get $|\tilde{a}_{bi}-a_{bi}|\leq \epsilon_{i_0,j_0}$.\\
			\item $a_{bi}\times a_{iw}> (a_{bi_0}+\epsilon_{i_0,j_0})\times(a_{i_0w}+\epsilon_{i_0,j_0})$\\\\
			In this case, the result can be proven by using a similar argument as in possibility (i).
		\end{enumerate}
		Hence the proof.
	\end{proof} 
\end{theorem}
\hspace{-0.75cm}
Theorem \ref{5exact_obj}, along with  the facts that $\epsilon_i=0$ for $i\in D_3$, $\epsilon_{i,j}=\epsilon_{j,i}$ for all $i,j\in D$, and that $\epsilon_{i,j}<\max\{\epsilon_i,\epsilon_j\}$ whenever $i,j\in D_1\cup D_3$ or $i,j\in D_2\cup D_3$, gives
\begin{equation}\label{5opt_obj_actual}
	\begin{split}
		\epsilon^*=\eta^*&=\max\{\epsilon_i,\epsilon_{i,j}:i,j\in D\}=\max\{\epsilon_i,\epsilon_{j,k}:i\in D_1\cup D_2, j\in D_1,k\in D_2\}\\
		&=\max\left\{\bigg|\frac{a_{bi}+a_{iw}+1-\sqrt{(a_{bi}+a_{iw}+1)^2-4(a_{bi}\times a_{iw}-a_{bw})}}{2}\bigg|,\right.\\
		&\quad\quad\quad\quad\quad\quad\quad\quad\quad\quad\quad\quad\quad\quad\quad\quad\quad\ \ \left.\bigg|\frac{a_{bi}\times a_{iw}-a_{bj}\times a_{jw}}{a_{bi}+a_{iw}+a_{bj}+a_{jw}}\bigg|:i,j\in D\right\},\\
		&=\max\left\{\bigg|\frac{a_{bi}+a_{iw}+1-\sqrt{(a_{bi}+a_{iw}+1)^2-4(a_{bi}\times a_{iw}-a_{bw})}}{2}\bigg|,\right.\\
		&\quad\quad\quad\quad\quad\quad\quad\quad\ \ \ \left.\frac{a_{bk}\times a_{kw}-a_{bj}\times a_{jw}}{a_{bj}+a_{jw}+a_{bk}+a_{kw}}:i\in D_1\cup D_2, j\in D_1,k\in D_2\right\},
	\end{split}
\end{equation}
which is analytical expression of $\epsilon^*$.
\begin{proposition}\label{5fix_values}
	Let $\epsilon_i$ and $\epsilon_{i,j}$ be as in equations \eqref{5CV} and \eqref{5CV'} respectively, let $\eta^*$ be the optimal objective value of problem \eqref{5optimization_2}, and let $(\tilde{A}_b^*,\tilde{A}_w^*)$ be an optimally modified PCS. Then 
\begin{equation}\label{5fix}
	\tilde{a}_{bw}^*=\begin{cases}
			a_{bw}-\epsilon_{i_0}\quad\quad\quad\quad\quad\quad\quad\quad\quad\quad\text{ if } \eta^*=\epsilon_{i_0}\text{ for some } i_0\in D_1,\\
			a_{bw}+\epsilon_{j_0}\quad\quad\quad\quad\quad\quad\quad\quad\quad\quad\text{ if } \eta^*=\epsilon_{j_0}\text{ for some } j_0\in D_2,\\
			(a_{bi_0}+\epsilon_{i_0,j_0})\times (a_{i_0w}+\epsilon_{i_0,j_0})\quad\ \text{ if } \eta^*=\epsilon_{i_0,j_0}\text{ for some } i_0,j_0\in D.
	\end{cases}
\end{equation}
\end{proposition}
\begin{proof}
	First, assume that $\eta^*=\epsilon_{i_0}\text{ for some } i_0\in D_1$. As discussed in Proposition \ref{5lower}, we get $0\leq\eta_{bi_0},\eta_{i_0w},\eta_{bw}\leq \eta^*=\epsilon_{i_0}$ such that $\tilde{a}_{bi_0}^*\in\{a_{bi_0}-\eta_{bi_0},a_{bi_0}+\eta_{bi_0}\}$, $\tilde{a}_{i_0w}^*\in\{a_{i_0w}-\eta_{i_0w},a_{i_0w}+\eta_{i_0w}\}$, and $\tilde{a}_{bw}^*\in\{a_{bw}-\eta_{bw},a_{bw}+\eta_{bw}\}$. Since $(\tilde{A}_b^*,\tilde{A}_w^*)$ is consistent, we have $\tilde{a}_{bi_0}^*\times \tilde{a}_{i_0w}^*=\tilde{a}_{bw}^*$. Now, it is easy to observe that if any one of $\eta_{bi_0}$, $\eta_{i_0w}$ or $\eta_{bw}$ is strictly less than $\epsilon_{i_0}$, then at least one of the remaining two must be strictly greater than $\epsilon_{i_0}$, which is not possible. This gives $\eta_{bi_0}=\eta_{i_0w}=\eta_{bw}=\epsilon_{i_0}$. It is clear that the only combination among the eight possible values of $(\tilde{a}_{bi_0}^*,\tilde{a}_{i_0w}^*,\tilde{a}_{bw}^*)$ that satisfies $\tilde{a}_{bi_0}^*\times\tilde{a}_{i_0w}^*=\tilde{a}_{bw}^*$ is $(a_{bi_0}+\epsilon_{i_0},a_{i_0w}+\epsilon_{i_0},a_{bw}-\epsilon_{i_0})$. Therefore, $\tilde{a}_{bw}^*=a_{bw}-\epsilon_{i_0}$. The result follows by an analogous argument if $\eta^*=\epsilon_{j_0}\text{ for some } j_0\in D_2$.\\\\
	Now, assume that $\eta^*=\epsilon_{i_0,j_0}\text{ for some } i_0,j_0\in D$. As discussed in Proposition \ref{5lower}, we get $0\leq\eta_{bi_0},\eta_{i_0w},\eta_{bj_0},\eta_{j_0w}\leq \eta^*=\epsilon_{i_0,j_0}$ such that $\tilde{a}_{bi_0}^*\in\{a_{bi_0}-\eta_{bi_0},a_{bi_0}+\eta_{bi_0}\}$, $\tilde{a}_{i_0w}^*\in\{a_{i_0w}-\eta_{i_0w},a_{i_0w}+\eta_{i_0w}\}$, $\tilde{a}_{bj_0}^*\in\{a_{bj_0}-\eta_{bj_0},a_{bj_0}+\eta_{bj_0}\}$, and $\tilde{a}_{j_0w}^*\in\{a_{j_0w}-\eta_{j_0w},a_{j_0w}+\eta_{j_0w}\}$. Since $(\tilde{A}_b^*,\tilde{A}_w^*)$ is consistent, we have $\tilde{a}_{bi_0}^*\times \tilde{a}_{i_0w}^*=\tilde{a}_{bj_0}^*\times \tilde{a}_{j_0w}^*$. Now, it is easy to observe that if any one of $\eta_{bi_0}$, $\eta_{i_0w}$, $\eta_{bj_0}$ or $\eta_{j_0w}$ is strictly less than $\epsilon_{i_0,j_0}$, then at least one of the remaining three must be strictly greater than $\epsilon_{i_0,j_0}$, which is not possible. This gives $\eta_{bi_0}=\eta_{i_0w}=\eta_{bj_0}=\eta_{j_0w}=\epsilon_{i_0,j_0}$. It is clear that the only combination among the sixteen possible values of $(\tilde{a}_{bi_0}^*,\tilde{a}_{i_0w}^*,\tilde{a}_{bj_0}^*,\tilde{a}_{j_0w}^*)$ that satisfies $\tilde{a}_{bi_0}^*\times\tilde{a}_{i_0w}^*=\tilde{a}_{bj_0}^*\times\tilde{a}_{j_0w}^*$ is $(a_{bi_0}+\epsilon_{i_0,j_0},a_{i_0w}+\epsilon_{i_0,j_0},a_{bj_0}-\epsilon_{i_0,j_0},a_{j_0w}-\epsilon_{i_0,j_0})$. Therefore, $\tilde{a}_{bw}^*=\tilde{a}_{bi_0}^*\times\tilde{a}_{i_0w}^*=(a_{bi_0}+\epsilon_{i_0,j_0})\times (a_{i_0w}+\epsilon_{i_0,j_0})$. Hence the proof.
\end{proof}
\hspace{-0.7cm}
Proposition \ref{5fix_values} implies that for any $(\tilde{A}_b^*,\tilde{A}_w^*)$, the value of $\tilde{a}_{bw}^*$ remains unchanged.\\\\
Using equations \eqref{5opt_obj_actual} and \eqref{5fix} in equation \eqref{5opt_weights}, we get optimal interval-weights.
\begin{theorem}\label{5best_pcs_actual}
	Let $\epsilon_i$ and $\epsilon_{i,j}$ be as in equations \eqref{5CV} and \eqref{5CV'} respectively, and let $\epsilon^*$ be the optimal objective value of problem \eqref{5optimization_1}. Then the following statements hold.
	\begin{enumerate}
		\item If $\epsilon^*=\epsilon_{i_0}$ for some $i_0\in D_1$, then $(\tilde{A}_b^*,\tilde{A}_w^*)$ given in equation \eqref{5best_pcs_1} is the only best optimally modified PCS.
		\item If $\epsilon^*=\epsilon_{j_0}$ for some $j_0\in D_2$, then $(\tilde{A}_b^*,\tilde{A}_w^*)$ given in equation \eqref{5best_pcs_2} is the only best optimally modified PCS.
		\item If $\epsilon^*=\epsilon_{i_0,j_0}$ for some $i_0,j_0\in D$, then $(\tilde{A}_b^*,\tilde{A}_w^*)$ given in equation \eqref{5best_pcs_3} is the only best optimally modified PCS.
	\end{enumerate}
\end{theorem}
\begin{proof}
	Let $(\tilde{A}_b'^*,\tilde{A}_w'^*)$ be an optimally modified PCS. Let $|\tilde{a}_{bi}'^*-a_{bi}|=\eta_{bi}$ and $|\tilde{a}_{iw}'^*-a_{iw}|=\eta_{iw}$, where $i\in D$. So, $0\leq\eta_{bi},\eta_{iw}\leq \eta^*$. Also, $\tilde{a}_{bi}'^*\in\{a_{bi}+\eta_{bi}, a_{bi}-\eta_{bi}\}$ and $\tilde{a}_{iw}'^*\in\{a_{iw}+\eta_{iw}, a_{iw}-\eta_{iw}\}$.\\\\
	First, assume that $\epsilon^*(=\eta^*)=\epsilon_{i_0}$ for some $i_0\in D_1$. By Proposition \ref{5fix_values}, $\tilde{a}_{bw}^*=\tilde{a}_{bw}'^*=a_{bw}-\epsilon_{i_0}$. Therefore, $|\tilde{a}_{bw}^*-a_{bw}|=|\tilde{a}_{bw}'^*-a_{bw}|$.\\\\
	Fix $i\in D$. If $a_{bi}\times a_{iw}\leq a_{bw}-\epsilon_{i_0}$, then by equation \eqref{5opt_dev_1}, $$|\tilde{a}_{bi}^*-a_{bi}|=|\tilde{a}_{iw}^*-a_{iw}|=\frac{-(a_{bi}+a_{iw})+\sqrt{(a_{bi}+a_{iw})^2-4\times a_{bi}\times a_{iw}+4(a_{bw}-\epsilon_{i_0})}}{2},$$ and thus, $$\max\{|\tilde{a}_{bi}^*-a_{bi}|,|\tilde{a}_{iw}^*-a_{iw}|\}=\frac{-(a_{bi}+a_{iw})+\sqrt{(a_{bi}+a_{iw})^2-4\times a_{bi}\times a_{iw}+4(a_{bw}-\epsilon_{i_0})}}{2}.$$ Since $(\tilde{A}_b'^*,\tilde{A}_w'^*)$ is consistent, we have $\tilde{a}_{bi}'\times \tilde{a}_{iw}'=\tilde{a}_{bw}'=a_{bw}-\epsilon_{i_0}$. This, along with \eqref{5dev_eq_1}, implies that if one of $\eta_{bi}$ and $\eta_{iw}$ is strictly less than $$\frac{-(a_{bi}+a_{iw})+\sqrt{(a_{bi}+a_{iw})^2-4\times a_{bi}\times a_{iw}+4(a_{bw}-\epsilon_{i_0})}}{2},$$ then the other necessarily exceeds it. This gives $\max\{|\tilde{a}_{bi}^*-a_{bi}|,|\tilde{a}_{iw}^*-a_{iw}|\}\leq \max\{\eta_{bi},\eta_{iw}\}=\max\{|\tilde{a}_{bi}'^*-a_{bi}|,|\tilde{a}_{iw}'^*-a_{iw}|\}$. Moreover, if $\max\{|\tilde{a}_{bi}^*-a_{bi}|,|\tilde{a}_{iw}^*-a_{iw}|\}=\max\{|\tilde{a}_{bi}'^*-a_{bi}|,|\tilde{a}_{iw}'^*-a_{iw}|\}$, then $$\eta_{bi}=\eta_{iw}=\frac{-(a_{bi}+a_{iw})+\sqrt{(a_{bi}+a_{iw})^2-4\times a_{bi}\times a_{iw}+4(a_{bw}-\epsilon_{i_0})}}{2},$$ which gives $\tilde{a}_{bi}'^*=\tilde{a}_{bi}^*$ and $\tilde{a}_{iw}'^*=\tilde{a}_{iw}^*$. Similar argument can be given if $a_{bi}\times a_{iw}> a_{bw}-\epsilon_{i_0}$. The conclusion also holds if $\eta^*=\epsilon_{j_0}\text{ for some } j_0\in D_2$ by an analogous reasoning.\\\\
	Now, assume that $\epsilon^*(=\eta^*)=\epsilon_{i_0,j_0}$ for some $i_0,j_0\in D$. By Proposition \ref{5fix_values}, $\tilde{a}_{bw}^*=\tilde{a}_{bw}'^*=(a_{bi_0}+\epsilon_{i_0,j_0})\times (a_{i_0w}+\epsilon_{i_0,j_0})$. So, $|\tilde{a}_{bw}^*-a_{bw}|=|\tilde{a}_{bw}'^*-a_{bw}|$.\\\\
	Fix $i\in D$. If $a_{bi}\times a_{iw}\leq (a_{bi_0}+\epsilon_{i_0,j_0})\times(a_{i_0w}+\epsilon_{i_0,j_0})$, then by equation \eqref{5opt_dev_3},
	\begin{equation*}
		\begin{split}
			|\tilde{a}_{bi}-a_{bi}|&=|\tilde{a}_{iw}-a_{iw}|\\
			&=\frac{-(a_{bi}+a_{iw})+\sqrt{(a_{bi}+a_{iw})^2-4\times a_{bi}\times a_{iw}+4(a_{bi_0}+\epsilon_{i_0,j_0})\times (a_{i_0w}+\epsilon_{i_0,j_0})}}{2},
		\end{split}
	\end{equation*}
	and thus, 
	\begin{equation*}
		\begin{split}
			\max&\{|\tilde{a}_{bi}-a_{bi}|,|\tilde{a}_{iw}-a_{iw}|\}\\
			&=\frac{-(a_{bi}+a_{iw})+\sqrt{(a_{bi}+a_{iw})^2-4\times a_{bi}\times a_{iw}+4(a_{bi_0}+\epsilon_{i_0,j_0})\times (a_{i_0w}+\epsilon_{i_0,j_0})}}{2}.
		\end{split}
	\end{equation*}
	Since $(\tilde{A}_b'^*,\tilde{A}_w'^*)$ is consistent, we have $\tilde{a}_{bi}'\times \tilde{a}_{iw}'=\tilde{a}_{bw}'=(a_{bi_0}+\epsilon_{i_0,j_0})\times (a_{i_0w}+\epsilon_{i_0,j_0})$. This, along with \eqref{5dev_eq_3}, implies that if one of $\eta_{bi}$ and $\eta_{iw}$ is strictly less than $$\frac{-(a_{bi}+a_{iw})+\sqrt{(a_{bi}+a_{iw})^2-4\times a_{bi}\times a_{iw}+4(a_{bi_0}+\epsilon_{i_0,j_0})\times (a_{i_0w}+\epsilon_{i_0,j_0})}}{2},$$ then the other necessarily exceeds it. This gives $\max\{|\tilde{a}_{bi}^*-a_{bi}|,|\tilde{a}_{iw}^*-a_{iw}|\}\leq \max\{\eta_{bi},\eta_{iw}\}=\max\{|\tilde{a}_{bi}'^*-a_{bi}|,|\tilde{a}_{iw}'^*-a_{iw}|\}$. Moreover, if $\max\{|\tilde{a}_{bi}^*-a_{bi}|,|\tilde{a}_{iw}^*-a_{iw}|\}=\max\{|\tilde{a}_{bi}'^*-a_{bi}|,|\tilde{a}_{iw}'^*-a_{iw}|\}$, then $$\eta_{bi}=\eta_{iw}=\frac{-(a_{bi}+a_{iw})+\sqrt{(a_{bi}+a_{iw})^2-4\times a_{bi}\times a_{iw}+4(a_{bi_0}+\epsilon_{i_0,j_0})\times (a_{i_0w}+\epsilon_{i_0,j_0})}}{2},$$ which gives $\tilde{a}_{bi}'^*=\tilde{a}_{bi}^*$ and $\tilde{a}_{iw}'^*=\tilde{a}_{iw}^*$. Similar argument can be given if $a_{bi}\times a_{iw}> (a_{bi_0}+\epsilon_{i_0,j_0})\times (a_{i_0w}+\epsilon_{i_0,j_0})$. Hence the proof.
\end{proof}
\hspace{-0.7cm}
Using this best optimally modified PCS, the best optimal weight set is obtained using equation \eqref{5equivalence}.
\subsubsection{Multiple Best or Worst Criteria}
In cases with multiple best/worst criteria, rather than arbitrarily designating one as the best or worst, we modify the original optimization model to assign equal weights to all best criteria and equal weights to all worst criteria.\\\\
Let $c_{b_1},c_{b_2},\ldots,c_{b_{n_1}}$ be the best and $c_{w_1},c_{w_2},\ldots,$ $c_{w_{n_2}}$ be the worst criteria. Therefore, $$a_{b_1i}=a_{b_2i}=\ldots=a_{b_{n_1}i}= a_{bi}\text{ (say)},\quad\quad a_{iw_1}=a_{iw_2}=\ldots=a_{iw_{n_2}}= a_{iw}\text{ (say)},$$ $$a_{b_1w_1}=\ldots =a_{b_1w_{n_2}}=a_{b_2w_1}=\ldots =a_{b_2w_{n_2}}=\ldots=a_{b_{n_1}w_1}=\ldots =a_{b_{n_1}w_{n_2}}= a_{bi}\text{ (say)}$$ for $i\in D'$, where $D'=C\setminus\{b_1,\ldots,b_{n_1},w_1,\ldots,w_{n_2}\}$.\\\\
To ensure $w_{b_1}=w_{b_2}=\ldots=w_{b_{n_1}}$ and $w_{w_1}=w_{w_2}=\ldots=w_{w_{n_2}}$, instead of considering the system of equations \eqref{5system}, consider the system of equations 
\begin{equation}\label{5system3}
	\begin{split}
		&\frac{w_{b_1}}{w_i}=\frac{w_{b_2}}{w_i}=\ldots=\frac{w_{b_{n_1}}}{w_i}=a_{bi},\quad \frac{w_i}{w_{w_1}}=\frac{w_i}{w_{w_2}}=\ldots=\frac{w_i}{w_{w_{n_2}}}=a_{iw},\\
		&\frac{w_{b_1}}{w_{w_1}}=\ldots=\frac{w_{b_1}}{w_{w_{n_2}}}=\frac{w_{b_2}}{w_{w_1}}=\ldots=\frac{w_{b_2}}{w_{w_{n_2}}}=\ldots=\frac{w_{b_{n_1}}}{w_{w_1}}=\ldots=\frac{w_{b_{n_1}}}{w_{w_{n_2}}}=a_{bw},\quad i\in D',\\
		&w_1+w_2+\ldots+w_n=1,
	\end{split}
\end{equation}
which is equivalent to the system of equations 
\begin{equation}\label{5system6}
	\begin{split}
		&\frac{w_{b_1}}{w_i}=a_{bi},\quad \frac{w_i}{w_{w_1}}=a_{iw},\quad \frac{w_{b_1}}{w_{w_1}}=a_{bw},\quad i\in D',\\
		&w_{b_1}=w_{b_2}=\ldots=w_{b_{n_1}},\quad w_{w_1}=w_{w_2}=\ldots=w_{w_{n_2}},\\
		&n_1 w_{b_1}+n_2 w_{w_1}+\sum_{k\in D'}w_k=1.
	\end{split}
\end{equation}
Consider the following minimization problem.
\begin{equation}\label{5optimization_7}
	\begin{split}
		&\min\epsilon \\
		&\text{subject to:}\\	
		&\left|\frac{w_{b_1}}{w_i}-a_{bi}  \right| \leq \epsilon, \quad \left|\frac{w_i}{w_{w_1}}-a_{iw}  \right| \leq \epsilon, \quad	\left| \frac{w_{b_1}}{w_{w_1}}-a_{bw} \right| \leq \epsilon, \\
		&n_1 w_{b_1}+n_2 w_{w_1}+\sum_{k\in D'}w_k=1,\quad w_j\geq 0\text{ for all } i\in D' \text{ and } j\in C.
	\end{split}
\end{equation}
Each optimal solution of problem \eqref{5optimization_7}, along with $w_{b_1}^*=w_{b_2}^*=\ldots=w_{b_{n_1}}^*$ and $w_{w_1}^*=w_{w_2}^*=\ldots=w_{w_{n_2}}^*$, gives an optimal weight set $W^*=\{w_1^*,w_2^*,\ldots,w_n^*\}$, and $\epsilon^*$ is a measurement of the accuracy of this weight set. Since $\epsilon^*$ is also the optimal objective value, it remains the same for all $W^*$.\\\\
When problem \eqref{5optimization_7} has multiple optimal solutions, the optimal interval-weight for criterion $c_k$ is $[w_k'^*,w_k''^*]$, where $w_k'^*$ and $w_k''^*$ represent the optimal objective values of problems
\begin{equation}\label{5optimization_8}
	\begin{split}
		&\min w_k \\
		&\text{subject to:}\\	
		&\left|\frac{w_{b_1}}{w_i}-a_{bi}  \right| \leq \epsilon^*, \quad \left|\frac{w_i}{w_{w_1}}-a_{iw}  \right| \leq \epsilon^*, \quad	\left| \frac{w_{b_1}}{w_{w_1}}-a_{bw} \right| \leq \epsilon^*, \\
		&n_1 w_{b_1}+n_2 w_{w_1}+\sum_{k\in D'}w_k=1,\quad w_j\geq 0\text{ for all } i\in D' \text{ and } j\in C
	\end{split}
\end{equation}
and
\begin{equation}\label{5optimization_9}
	\begin{split}
		&\max w_k \\
		&\text{subject to:}\\	
		&\left|\frac{w_{b_1}}{w_i}-a_{bi}  \right| \leq \epsilon^*, \quad \left|\frac{w_i}{w_{w_1}}-a_{iw}  \right| \leq \epsilon^*, \quad	\left| \frac{w_{b_1}}{w_{w_1}}-a_{bw} \right| \leq \epsilon^*, \\
		&n_1 w_{b_1}+n_2 w_{w_1}+\sum_{k\in D'}w_k=1,\quad w_j\geq 0\text{ for all } i\in D' \text{ and } j\in C
	\end{split}
\end{equation}
respectively.\\\\
Consider the following minimization problem.
\begin{equation}\label{5optimization_10}
	\begin{split}
		&\min\eta \\
		&\text{subject to:}\\	
		&\left|\tilde{a}_{bi}-a_{bi}  \right| \leq \eta, \quad \left|\tilde{a}_{iw}-a_{iw}  \right| \leq \eta, \quad	\left| \tilde{a}_{bw}-a_{bw} \right| \leq \eta, \\
		&\tilde{a}_{bi}\times\tilde{a}_{iw}=\tilde{a}_{bw},\quad \tilde{a}_{bi},\tilde{a}_{iw},\tilde{a}_{bw}\geq 0\text{ for all } i\in D'.
	\end{split}
\end{equation}
Each optimal solution of problem \eqref{5optimization_10}, along with $\tilde{a}_{bb}^*=\tilde{a}_{ww}^*=1$, leads to an optimally modified PCS $(\tilde{A}_b^*,\tilde{A}_w^*)$ and $\eta^*$ is a measurement of its accuracy. Since $\eta^*$ is also the optimal objective value, it remains the same for all $(\tilde{A}_b^*,\tilde{A}_w^*)$.\\\\
Analogous to the case with unique best and worst criteria, it can be verified that equation \eqref{5equivalence} holds for the collections of optimal solutions of problems \eqref{5optimization_7} and \eqref{5optimization_10}, and thus, the two problems are equivalent.
\begin{theorem}\label{5op_modified_PCS_multi}
	Let $\epsilon_i$ and $\epsilon_{i,j}$ be as in equations \eqref{5CV} and \eqref{5CV'} respectively, and let $\epsilon^*$ and $\eta^*$ be the optimal objective values of problems \eqref{5optimization_7} and \eqref{5optimization_10}. Then
	\begin{equation}\label{5opt_obj_multi}
		\epsilon^*=\eta^*=\max\{\epsilon_i,\epsilon_{i,j}:i,j\in D\}.
	\end{equation}
Also, the collection of all optimally modified PCS is
\begin{equation}\label{5op_PCS_multi}
	\begin{split}
		&\tilde{a}_{bw}^*=\begin{cases}
			a_{bw}-\epsilon_{i_0}\quad\quad\quad\quad\quad\quad\quad\quad\quad\quad\text{ if } \eta^*=\epsilon_{i_0}\text{ for some } i_0\in D_1,\\
			a_{bw}+\epsilon_{j_0}\quad\quad\quad\quad\quad\quad\quad\quad\quad\quad\text{ if } \eta^*=\epsilon_{j_0}\text{ for some } j_0\in D_2,\\
			(a_{bi_0}+\epsilon_{i_0,j_0})\times (a_{i_0w}+\epsilon_{i_0,j_0})\quad\ \text{ if } \eta^*=\epsilon_{i_0,j_0}\text{ for some } i_0,j_0\in D.
		\end{cases}\\
		&\tilde{a}_{iw}^*\in \bigg[\max\biggl\{a_{iw}-\epsilon^*,\frac{\tilde{a}_{bw}^*}{a_{bi}+\epsilon^*}\biggr\},\min\biggl\{a_{iw}+\epsilon^*,\frac{\tilde{a}_{bw}^*}{a_{bi}-\epsilon^*}\biggr\}\bigg]\text{ with}\\
		&\tilde{a}_{bi}^*=\frac{\tilde{a}_{bw}^*}{\tilde{a}_{iw}^*}\text{ for all } i\in D'.
	\end{split}
\end{equation}
\end{theorem}
\begin{proof}
	The result can be proven by replicating the arguments of Theorem \ref{5exact_obj} and Proposition \ref{5fix_values}, and is thus omitted.
\end{proof}
\hspace{-0.7cm}
Using equation \eqref{5op_PCS_multi} in equation \eqref{5opt_weights_1}, we get
\begin{equation}\label{5opt_weights_multi}
	\begin{split}
		{w_{b_1}^l}^*&={w_{b_2}^l}^*=\ldots={w_{b_{n_1}}^l}^*=\frac{\tilde{a}_{bw}^*}{n_2+n_1\tilde{a}_{bw}^*+\displaystyle\sum_{j\in D'}\min\biggl\{a_{jw}+\epsilon^*,\frac{\tilde{a}_{bw}^*}{a_{bj}-\epsilon^*}\biggr\}},\\
		{w_{b_1}^u}^*&={w_{b_2}^u}^*=\ldots={w_{b_{n_1}}^u}^*=\frac{\tilde{a}_{bw}^*}{n_2+n_1\tilde{a}_{bw}^*+\displaystyle\sum_{j\in D'}\max\biggl\{a_{jw}-\epsilon^*,\frac{\tilde{a}_{bw}^*}{a_{bj}+\epsilon^*}\biggr\}},\\
		{w_{w_1}^l}^*&={w_{w_2}^l}^*=\ldots={w_{w_{n_2}}^l}^*=\frac{1}{n_2+n_1\tilde{a}_{bw}^*+\displaystyle\sum_{j\in D'}\min\biggl\{a_{jw}+\epsilon^*,\frac{\tilde{a}_{bw}^*}{a_{bj}-\epsilon^*}\biggr\}},\\
		{w_{w_1}^u}^*&={w_{w_2}^u}^*=\ldots={w_{w_{n_2}}^u}^*=\frac{1}{n_2+n_1\tilde{a}_{bw}^*+\displaystyle\sum_{j\in D'}\max\biggl\{a_{jw}-\epsilon^*,\frac{\tilde{a}_{bw}^*}{a_{bj}+\epsilon^*}\biggr\}},\\
		{w_i^l}^*&=\frac{\max\biggl\{a_{iw}-\epsilon^*,\frac{\tilde{a}_{bw}^*}{a_{bi}+\epsilon^*}\biggr\}}{n_2+n_1\tilde{a}_{bw}^*+\max\biggl\{a_{iw}-\epsilon^*,\frac{\tilde{a}_{bw}^*}{a_{bi}+\epsilon^*}\biggr\}+\displaystyle\sum_{\substack{j\in D'\\j\neq i}}\min\biggl\{a_{jw}+\epsilon^*,\frac{\tilde{a}_{bw}^*}{a_{bj}-\epsilon^*}\biggr\}},\\
		{w_i^u}^*&=\frac{\min\biggl\{a_{iw}+\epsilon^*,\frac{\tilde{a}_{bw}^*}{a_{bi}-\epsilon^*}\biggr\}}{n_2+n_1\tilde{a}_{bw}^*+\min\biggl\{a_{iw}+\epsilon^*,\frac{\tilde{a}_{bw}^*}{a_{bi}-\epsilon^*}\biggr\}+\displaystyle\sum_{\substack{j\in D'\\j\neq i}}\max\biggl\{a_{jw}-\epsilon^*,\frac{\tilde{a}_{bw}^*}{a_{bj}+\epsilon^*}\biggr\}},\quad \text{where } i\in D'.
	\end{split}
\end{equation}
Equations \eqref{5opt_obj_actual} and \eqref{5opt_obj_multi} imply that for a PCS involving multiple best/worst criteria, both optimization models \eqref{5optimization_1} and \eqref{5optimization_7} have the same optimal objective values. This indicates that optimization model \eqref{5optimization_7} selects the optimal weight sets in which all best criteria share equal weights and all worst criteria share equal weights.
\begin{theorem}\label{5exact_obj_1}
	Let $\epsilon_i$ and $\epsilon_{i,j}$ be as in equations \eqref{5CV} and \eqref{5CV'} respectively, and let $\epsilon^*$ be the optimal objective value of problem \eqref{5optimization_7}. Then the following statements hold.
	\begin{enumerate}
		\item If $\epsilon^*=\epsilon_{i_0}$ for some $i_0\in D_1$, then $(\tilde{A}_b,\tilde{A}_w)$ defined as 
		\begin{equation}\label{5best_pcs_4}
			\begin{split}
				&\tilde{a}_{bi}=\frac{a_{bi}-a_{iw}+\sqrt{(a_{bi}+a_{iw})^2-4\times a_{bi}\times a_{iw}+4(a_{bw}-\epsilon_{i_0})}}{2},\\
				&\tilde{a}_{iw}=\frac{-a_{bi}+a_{iw}+\sqrt{(a_{bi}+a_{iw})^2-4\times a_{bi}\times a_{iw}+4(a_{bw}-\epsilon_{i_0})}}{2}\quad\text{for all } i\in D',\\
				&\tilde{a}_{bw}=a_{bw}-\epsilon_{i_0}
			\end{split}
		\end{equation}
		is the only best optimally modified PCS.
		\item If $\epsilon^*=\epsilon_{j_0}$ for some $j_0\in D_2$, then $(\tilde{A}_b,\tilde{A}_w)$ defined as 
		\begin{equation}\label{5best_pcs_5}
			\begin{split}
				&\tilde{a}_{bi}=\frac{a_{bi}-a_{iw}+\sqrt{(a_{bi}+a_{iw})^2-4\times a_{bi}\times a_{iw}+4(a_{bw}+\epsilon_{j_0})}}{2},\\
				&\tilde{a}_{iw}=\frac{-a_{bi}+a_{iw}+\sqrt{(a_{bi}+a_{iw})^2-4\times a_{bi}\times a_{iw}+4(a_{bw}+\epsilon_{j_0})}}{2}\quad\text{for all } i\in D',\\
				&\tilde{a}_{bw}=a_{bw}+\epsilon_{j_0}
			\end{split}
		\end{equation}
		is the only best optimally modified PCS.
		\item If $\epsilon^*=\epsilon_{i_0,j_0}$ for some $i_0,j_0\in D$, then $(\tilde{A}_b,\tilde{A}_w)$ defined as 
		\begin{equation}\label{5best_pcs_6}
			\begin{split}
				&\tilde{a}_{bi}=\frac{a_{bi}-a_{iw}+\sqrt{(a_{bi}+a_{iw})^2-4\times a_{bi}\times a_{iw}+4(a_{bi_0}+\epsilon_{i_0,j_0})\times (a_{i_0w}+\epsilon_{i_0,j_0})}}{2},\\
				&\tilde{a}_{iw}=\frac{-a_{bi}+a_{iw}+\sqrt{(a_{bi}+a_{iw})^2-4\times a_{bi}\times a_{iw}+4(a_{bi_0}+\epsilon_{i_0,j_0})\times (a_{i_0w}+\epsilon_{i_0,j_0})}}{2}\\
				&\quad\quad\quad\quad\quad\quad\quad\quad\quad\quad\quad\quad\quad\quad\quad\quad\quad\quad\quad\quad\quad\quad\quad\quad\quad\quad\quad\quad\quad\text{for all } i\in D',\\
				&\tilde{a}_{bw}=(a_{bi_0}+\epsilon_{i_0,j_0})\times (a_{i_0w}+\epsilon_{i_0,j_0})
			\end{split}
		\end{equation}
		is the only best optimally modified PCS.
	\end{enumerate}
	\end{theorem}
	\hspace{-0.725cm}
	Using the best optimally modified PCS, the best optimal weight set is determined using equation \eqref{5equivalence}.
\subsection{Consistency Analysis}
\begin{figure}[H]
	\centering
	\includegraphics[height=10cm,width=15cm]{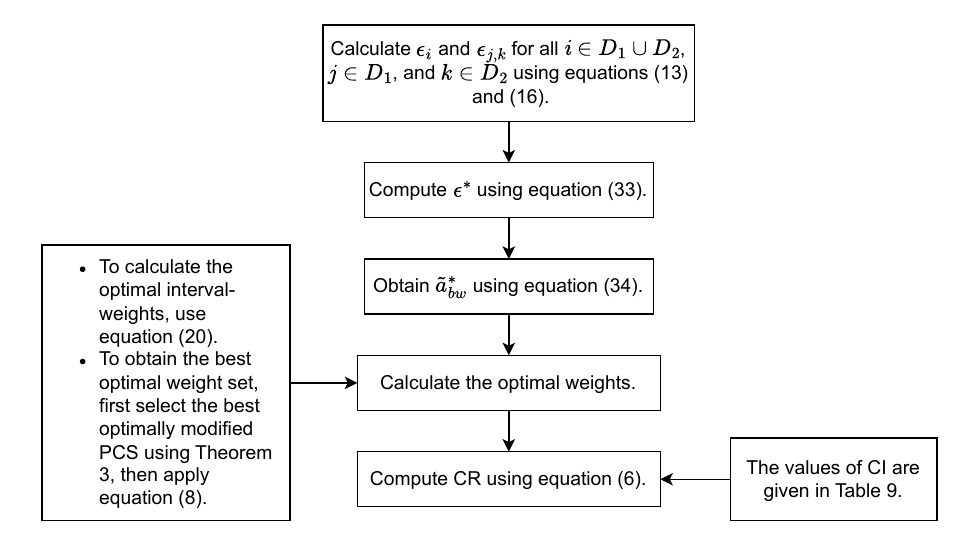}
	\caption{Workflow for weight computation and consistency analysis (unique best and worst criterion)}
	\label{5flowchart_1}		
\end{figure}

\begin{figure}[H]
	\centering
	\includegraphics[height=10cm,width=15cm]{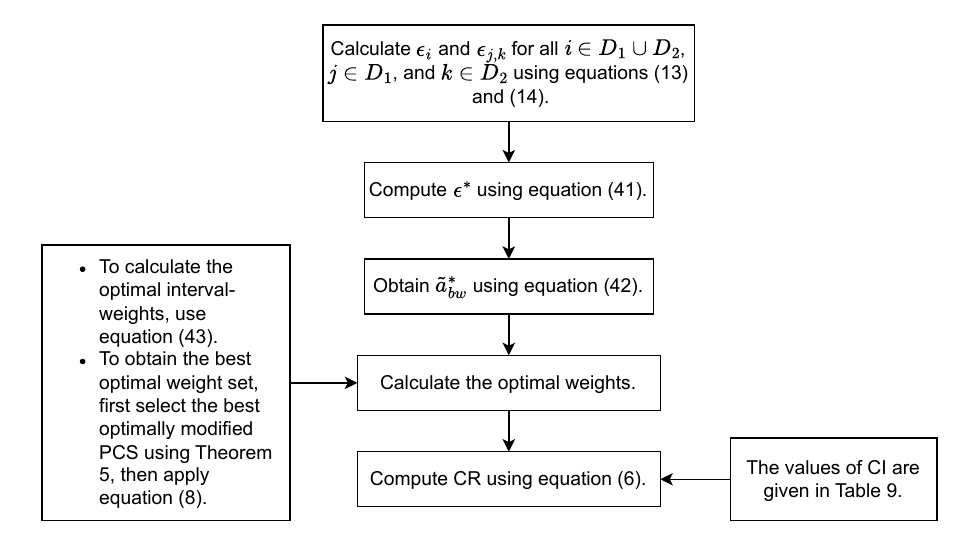}
	\caption{Workflow for weight computation and consistency analysis (multiple best/worst criteria)}
	\label{5flowchart_2}		
\end{figure}

This subsection derives analytical expressions for the CI and CR.
\begin{proposition}\label{5upper_1}
	Let $i,j\in D_1$ be such that $a_{bi}\leq a_{bj}$ and $a_{iw}\leq a_{jw}$. Then $\epsilon_j \leq \epsilon_i$.	
\end{proposition}
\begin{proof}
	Suppose, if possible, $\epsilon_i<\epsilon_j$. Then we get $a_{bi}+\epsilon_i<a_{bj}+\epsilon_j$ and $a_{iw}+\epsilon_i<a_{jw}+\epsilon_j$. This gives
	\begin{equation}\label{5inequality}
		(a_{bi}+\epsilon_i)\times(a_{iw}+\epsilon_i)<(a_{bj}+\epsilon_j)\times(a_{jw}+\epsilon_j).
	\end{equation} 
	From equation \eqref{5CV1}, we have $(a_{bi}+\epsilon_i)\times(a_{iw}+\epsilon_i)=a_{bw}-\epsilon_i$ and $(a_{bj}+\epsilon_j)\times(a_{jw}+\epsilon_j)=a_{bw}-\epsilon_j$. Therefore, by equation \eqref{5inequality}, we get $a_{bw}-\epsilon_i<a_{bw}-\epsilon_j$, which is contradiction as $\epsilon_i<\epsilon_j$ implies $a_{bw}-\epsilon_j<a_{bw}-\epsilon_i$. Thus, $\epsilon_j \leq \epsilon_i$. Hence the proof.
\end{proof}
\hspace{-0.725cm}
Proposition \ref{5upper_1}, along with the fact that $a_{bi},a_{iw}\geq 1$ for all $i\in D_1$, implies that for $(A_b,A_w)$ with $n\geq 3$,
\begin{equation}\label{5cv1}
	\max\{\epsilon_i:i\in D_1\}=\frac{-3+\sqrt{4a_{bw}+5}}{2}.
\end{equation}
\begin{proposition}\label{5upper_2}
	Let $i,j\in D_2$ be such that $a_{bi}\leq a_{bj}$ and $a_{iw}\leq a_{jw}$. Then $\epsilon_i \leq \epsilon_j$.	
\end{proposition}
\begin{proof}
	The proof is similar to the proof of Proposition \ref{5upper_1}, and thus omitted.
\end{proof}
\hspace{-0.725cm}
Proposition \ref{5upper_2}, along with the fact that $1\leq a_{bi},a_{iw}\leq a_{bw}$ for all $i\in D_2$, implies that for $(A_b,A_w)$ with $n\geq 3$,
\begin{equation}\label{5cv2}
	\max\{\epsilon_i:i\in D_2\}=\frac{2a_{bw}+1-\sqrt{8a_{bw}+1}}{2}.
\end{equation}
\begin{proposition}\label{5upper_3}
	Let $i,j,k\in D$ be such that $a_{bi}\leq a_{bj}\leq a_{bk}$ and $a_{iw}\leq a_{jw}\leq a_{kw}$. Then $\epsilon_{i,j} \leq \epsilon_{i,k}$ and $\epsilon_{j,k}\leq \epsilon_{i,k}$.	
\end{proposition}
\begin{proof}
	The proof is similar to the proof of Proposition \ref{5upper_1}, and thus omitted.
\end{proof}
\hspace{-0.725cm}
Proposition \ref{5upper_3}, along with the fact that $1\leq a_{bi},a_{iw}\leq a_{bw}$ for all $i\in D$, implies that for $(A_b,A_w)$ with $n\geq 4$,
\begin{equation}\label{5cv3}
	\max\{\epsilon_{j,k}:j\in D_1,k\in D_2\}=\frac{a_{bw}^2-1}{2a_{bw}+2}.
\end{equation}
Using equations \eqref{5cv1}, \eqref{5cv2}, and \eqref{5cv3} in equation \eqref{5opt_obj_actual}, we get $$\text{CI}_{a_{bw}}(n)=\begin{cases}
	\max\biggl\{\frac{-3+\sqrt{4a_{bw}+5}}{2},\frac{2a_{bw}+1-\sqrt{8a_{bw}+1}}{2}\biggr\}\quad\quad\quad\quad\ \text{if } n=3,\\
	\max\biggl\{\frac{-3+\sqrt{4a_{bw}+5}}{2},\frac{2a_{bw}+1-\sqrt{8a_{bw}+1}}{2},\frac{a_{bw}^2-1}{2a_{bw}+2}\biggr\}\quad \text{if } n\geq4.
\end{cases}$$
Let $$f(x)=\frac{2a_{bw}+1-\sqrt{8a_{bw}+1}}{2}-\frac{-3+\sqrt{4a_{bw}+5}}{2},\quad x\in [1,\infty).$$ Note that $$f'(x)=1-\frac{2}{\sqrt{8a_{bw}+1}}-\frac{1}{\sqrt{4a_{bw}+5}}.$$ Since $a_{bw}\geq 1$, we get $f'(x)\geq 0$, i.e., $f$ is increasing. Therefore, $f(x)\geq f(1)=0$ for all $x\in [1,\infty).$ This gives 
\begin{equation}\label{5ci_reduce}
	\frac{2a_{bw}+1-\sqrt{8a_{bw}+1}}{2}\geq\frac{-3+\sqrt{4a_{bw}+5}}{2},
\end{equation} and thus, 
\begin{equation}\label{5ci_ana}
	\text{CI}_{a_{bw}}(n)=\begin{cases}
		\frac{2a_{bw}+1-\sqrt{8a_{bw}+1}}{2}\quad \quad\quad\  \quad\quad\quad\quad\ \text{if } n=3,\\
		\max\biggl\{\frac{2a_{bw}+1-\sqrt{8a_{bw}+1}}{2},\frac{a_{bw}^2-1}{2a_{bw}+2}\biggr\}\quad \text{if } n\geq4.
	\end{cases}
\end{equation}
equation \eqref{5ci_ana} is an analytical expression for the CI. Table \ref{5ci_values} presents the values of CI for the scales given in Table \ref{5scale}.\\
\begin{table}[H]
	\caption {The values of CI for the scales listed in Table \ref{5scale} \label{5ci_values}}
	\centering
	\scalebox{0.95}{\footnotesize
		\begin{tabular}{ccccccccccccccc}
			\toprule[0.1em]
			\multicolumn{3}{c}{Saaty}&\phantom{}&\multicolumn{3}{c}{Salo-H{\"a}m{\"a}l{\"a}inen}&\phantom{}&\multicolumn{3}{c}{Lootsma}&\phantom{}&\multicolumn{3}{c}{Donegan-Dodd-}\\
			\multicolumn{3}{c}{scale}&&\multicolumn{3}{c}{scale}&&\multicolumn{3}{c}{scale}&&\multicolumn{3}{c}{McMaster scale}\\
			\cline{1-3}\cline{5-7}\cline{9-11}\cline{13-15}
			$a_{bw}$&\multicolumn{2}{c}{CI}&&$a_{bw}$&\multicolumn{2}{c}{CI}&&$a_{bw}$&\multicolumn{2}{c}{CI}&&$a_{bw}$&\multicolumn{2}{c}{CI}\\
			\cline{2-3}\cline{6-7}\cline{10-11}\cline{14-15}
			&$n=3$&$n\geq 4$&&&$n=3$&$n\geq 4$&&&$n=3$&$n\geq 4$&&&$n=3$&$n\geq 4$\\
			\midrule
			$2$&$0.4384$&$0.5$&&$1.2222$&$0.0807$&$0.1111$&&$\sqrt{2}$&$0.1597$&$0.2071$&&$1.1257$&$0.0441$&$0.0629$\\
			$3$&$1$&$1$&&$1.5$&$0.1972$&$0.25$&&$2$&$0.4384$&$0.5$&&$1.2715$&$0.1003$&$0.1358$\\
			$4$&$1.6277$&$1.6277$&&$1.8571$&$0.3661$&$0.4286$&&$2\sqrt{2}$&$0.8980$&$0.9142$&&$1.4470$&$0.1739$&$0.2235$\\
			$5$&$2.2984$&$2.2984$&&$2.3333$&$0.6160$&$0.6667$&&$4$&$1.6277$&$1.6277$&&$1.6684$&$0.2745$&$0.3342$\\
			$6$&$3$&$3$&&$3$&$1$&$1$&&$4\sqrt{2}$&$2.7563$&$2.7563$&&$1.9670$&$0.4215$&$0.4835$\\
			$7$&$3.7250$&$3.7250$&&$4$&$1.6277$&$1.6277$&&$8$&$4.4688$&$4.4688$&&$2.4142$&$0.6607$&$0.7071$\\
			$8$&$4.4688$&$4.4688$&&$5.6667$&$2.7633$&$2.7633$&&$8\sqrt{2}$&$7.0307$&$7.0307$&&$3.2289$&$1.1390$&$1.1390$\\
			$9$&$5.2279$&$5.2279$&&$9$&$5.2279$&$5.2279$&&$16$&$10.8211$&$10.8211$&&$5.8284$&$2.8778$&$2.8778$\\
			\bottomrule[0.1em]
	\end{tabular}}
\end{table}
\hspace{-0.725cm}
By substituting the expressions for $\epsilon^*$ from equation \eqref{5opt_obj_actual} and CI from equation \eqref{5ci_ana} into equation \eqref{5CR}, we get the analytical expression for CR.\\\\
Some key points regarding the CI and CR are as follows.
\begin{enumerate}
	\item For the Saaty scale, we have CI$_2(4)=0.5$, which is compatible with Example 5 discussed in Section 3.
	\item Since optimization models \eqref{5optimization_1} and \eqref{5optimization_7} produce the same optimal objective values for any PCS, equation \eqref{5ci_ana} remains valid even when multiple best/worst criteria are present.
	\item For a consistency indicator to demonstrate reasonable behavior, it must satisfy some specific properties \cite{koczkodaj2018axiomatization,koczkodaj2014axiomatization,brunelli2024inconsistency}. Proposition \ref{5cr_pro} outlines several of these properties without proof, as their proofs are analogous to those in \cite[Proposition 1]{liang2020consistency}.
\end{enumerate}
\begin{proposition}\label{5cr_pro}
	CR exhibits the following properties.
	\begin{enumerate}
		\item CR is normalized, i.e., $0\leq$ CR $\leq 1$.
		\item CR $=0$ if and only if $(A_b,A_w)$ is consistent.
		\item CR exhibits permutation invariance with respect to the criteria indices.
		\item CR is non-increasing with respect to criterion elimination.
		\item CR is a continuous function of $a_{bi}$, $a_{iw}$, and $a_{bw}$.
		\item For a consistent $(A_b,A_w)$, CR increases when either $a_{bi}$ or $a_{iw}$ moves away from its original value in the range $[1,a_{bw}]$.
	\end{enumerate}
\end{proposition}
Fig. \ref{5flowchart_1} shows the process for calculating optimal weights and analyzing consistency when a unique best and worst criterion exist, whereas Fig. \ref{5flowchart_2} shows the corresponding process in the presence of multiple best/worst criteria.
\subsection{Numerical Examples}
In this subsection, we revisit the six examples (Example 1–Example 6) from the section ``Research Gap" to demonstrate and validate the proposed approach.\\\\
\textbf{Example 1:} Here, $D_1=\{4\}$, $D_2=\{2,3\}$, and $D_3=\emptyset$. By equations \eqref{5CV} and \eqref{5CV'}, we get $\epsilon_2=0.4056$, $\epsilon_3=0.3944$, $\epsilon_4=0.5363$, $\epsilon_{4,2}=0.5163$, and $\epsilon_{4,3}=0.5422$. So, from equation \eqref{5opt_obj_actual}, it follows that $\epsilon^*=\max\{\epsilon_2,\epsilon_3,\epsilon_4,\epsilon_{4,2},\epsilon_{4,3}\}=\epsilon_{4,3}=0.5422$. Now, equation \ref{5fix} gives $\tilde{a}_{15}^*=(a_{14}+\epsilon_{4,3})\times(a_{45}+\epsilon_{4,3})=8.4987$. From equation \eqref{5opt_weights}, we get
\begin{equation*}
	\begin{split}
	{w_1^l}^*&=\frac{\tilde{a}_{15}^*}{1+\tilde{a}_{15}^*+\displaystyle\sum_{j=2,3,4}\min\biggl\{a_{j5}+\epsilon^*,\frac{\tilde{a}_{15}^*}{a_{1j}-\epsilon^*}\biggr\}}\\
	&=\frac{8.4987}{1+8.4987+1.0048+3.4578+3.5422}\\
	&=0.4855.
	\end{split}
\end{equation*}
Similarly, ${w_1^u}^*=0.4868$. Therefore, $w_1^*=[0.4855,0.4868]$. Computing in the same fashion, we obtain $w_2^*=[0.0549,0.0574]$, $w_3^*=[0.1975,0.1981]$, $w_4^*=[0.2024,0.2029]$, and $w_5^*=[0.0571,0.0573]$.\\\\
Since $\epsilon^*=\epsilon_{4,3}$, by statement \textit{3} of Theorem \ref{5best_pcs_actual}, 
\begin{equation*}
	\begin{split}
		&\tilde{a}_{bi}^*=\frac{a_{1i}-a_{i5}+\sqrt{(a_{1i}+a_{i5})^2-4\times a_{1i}\times a_{i5}+4(a_{14}+\epsilon_{4,3})\times (a_{45}+\epsilon_{4,3})}}{2},\\
		&\tilde{a}_{iw}^*=\frac{-a_{1i}+a_{i5}+\sqrt{(a_{1i}+a_{i5})^2-4\times a_{1i}\times a_{i5}+4(a_{14}+\epsilon_{4,3})\times (a_{45}+\epsilon_{4,3})}}{2},
	\end{split}
\end{equation*}
where $i=2,3,4$, along with $\tilde{a}_{11}^*=\tilde{a}_{55}^*=1$ and $\tilde{a}_{15}^*=8.4987$, form the best optimally modified PCS. Therefore, $$\tilde{A}_b^*=(1,8.4962,2.4519,2.3934,8.4638),\quad \tilde{A}_w^*=(8.4638,0.9962,3.4519,3.5363,1)^T$$ is the best optimally modified PCS. Thus, by equation \eqref{5equivalence}, $\{0.4857,0.0571,0.1976,0.2024,0.0572\}$ is the best optimal weight set. Now, by equation \eqref{5CR}, we get CR $=0.1037$.\\\\
Table \ref{5example_1} shows that the obtained $\epsilon^*$, interval-weights, and the best optimal weight set coincide with their actual values, validating the proposed framework.\\\\
\textbf{Example 2:} Here, $D_1=\{4\}$, $D_2=\{2,3\}$, and $D_3=\emptyset$. By equations \eqref{5CV} and \eqref{5CV'}, we get $\epsilon_2=1.6054$, $\epsilon_3=1.4764$, $\epsilon_4=1.7228$, $\epsilon_{4,2}=1.7882$, and $\epsilon_{4,3}=1.7999$. So, from equation \eqref{5opt_obj_actual}, it follows that $\epsilon^*=\max\{\epsilon_2,\epsilon_3,\epsilon_4,\epsilon_{4,2},\epsilon_{4,3}\}=\epsilon_{4,3}=1.7999$. Now, equation \ref{5fix} gives $\tilde{a}_{15}^*=(a_{14}+\epsilon_{4,3})\times(a_{45}+\epsilon_{4,3})=14.8760$. From equation \eqref{5opt_weights}, we get
\begin{equation*}
	\begin{split}
		{w_1^l}^*&=\frac{\tilde{a}_{15}^*}{1+\tilde{a}_{15}^*+\displaystyle\sum_{j=2,3,4}\min\biggl\{a_{j5}+\epsilon^*,\frac{\tilde{a}_{15}^*}{a_{1j}-\epsilon^*}\biggr\}}\\
		&=\frac{14.8760}{1+14.8760+1.0476+3.8569+3.2141}\\
		&=0.6200.
	\end{split}
\end{equation*}
Similarly, ${w_1^u}^*=0.6205$. Therefore, $w_1^*=[0.6200,0.6205]$. Computing in the same fashion, we obtain $w_2^*=[0.0429,0.0437]$, $w_3^*=[0.1607,0.1609]$, $w_4^*=[0.1340,0.1341]$, and $w_5^*=[0.0416,0.0417]$.\\\\
Since $\epsilon^*=\epsilon_{4,3}$, by statement \textit{3} of Theorem \ref{5best_pcs_actual}, 
\begin{equation*}
	\begin{split}
		&\tilde{a}_{bi}^*=\frac{a_{1i}-a_{i5}+\sqrt{(a_{1i}+a_{i5})^2-4\times a_{1i}\times a_{i5}+4(a_{14}+\epsilon_{4,3})\times (a_{45}+\epsilon_{4,3})}}{2},\\
		&\tilde{a}_{iw}^*=\frac{-a_{1i}+a_{i5}+\sqrt{(a_{1i}+a_{i5})^2-4\times a_{1i}\times a_{i5}+4(a_{14}+\epsilon_{4,3})\times (a_{45}+\epsilon_{4,3})}}{2},
	\end{split}
\end{equation*}
where $i=2,3,4$, along with $\tilde{a}_{11}^*=\tilde{a}_{55}^*=1$ and $\tilde{a}_{15}^*=14.8760$, form the best optimally modified PCS. Therefore, $$\tilde{A}_b^*=(1,14.2179,3.8569,4.6283,14.8760), \quad \tilde{A}_w^*=(14.8760,1.0463,3.8569,3.2141,1)^T$$ is the best optimally modified PCS. Thus, by equation \eqref{5equivalence}, $\{0.6200,0.0436,0.1607,0.1340,0.0417\}$ is the best optimal weight set. Now, by equation \eqref{5CR}, we get CR $=0.1663$.\\\\
Table \ref{5example_2} shows that the obtained $\epsilon^*$, interval-weights, and the best optimal weight set coincide with their actual values, validating the proposed framework.\\\\
\textbf{Example 3:} Here, $D_1=\{4\}$, $D_2=\{2,3\}$, and $D_3=\emptyset$. By equations \eqref{5CV} and \eqref{5CV'}, we get $\epsilon_2=0.6959$, $\epsilon_3=0.6781$, $\epsilon_4=0.8086$, $\epsilon_{4,2}=0.8825$, and $\epsilon_{4,3}=0.8975$. So, from equation \eqref{5opt_obj_actual}, it follows that $\epsilon^*=\max\{\epsilon_2,\epsilon_3,\epsilon_4,\epsilon_{4,2},\epsilon_{4,3}\}=\epsilon_{4,3}=0.8975$. Now, equation \ref{5fix} gives $\tilde{a}_{15}^*=(a_{14}+\epsilon_{4,3})\times(a_{45}+\epsilon_{4,3})=5.4354$. From equation \eqref{5opt_weights}, we get
\begin{equation*}
	\begin{split}
		{w_1^l}^*&=\frac{\tilde{a}_{15}^*}{1+\tilde{a}_{15}^*+\displaystyle\sum_{j=2,3,4}\min\biggl\{a_{j5}+\epsilon^*,\frac{\tilde{a}_{15}^*}{a_{1j}-\epsilon^*}\biggr\}}\\
		&=\frac{5.4354}{1+5.4354+1.1023+2.3314+2.8645}\\
		&=0.4269.
	\end{split}
\end{equation*}
Similarly, ${w_1^u}^*=0.4280$. Therefore, $w_1^*=[0.4269,0.4280]$. Computing in the same fashion, we obtain $w_2^*=[0.0842,0.0866]$, $w_3^*=[0.1831,0.1836]$, $w_4^*=[0.2250,0.2255]$, and $w_5^*=[0.0785,0.0787]$.\\\\
Since $\epsilon^*=\epsilon_{4,3}$, by statement \textit{3} of Theorem \ref{5best_pcs_actual}, 
\begin{equation*}
	\begin{split}
		&\tilde{a}_{bi}^*=\frac{a_{1i}-a_{i5}+\sqrt{(a_{1i}+a_{i5})^2-4\times a_{1i}\times a_{i5}+4(a_{14}+\epsilon_{4,3})\times (a_{45}+\epsilon_{4,3})}}{2},\\
		&\tilde{a}_{iw}^*=\frac{-a_{1i}+a_{i5}+\sqrt{(a_{1i}+a_{i5})^2-4\times a_{1i}\times a_{i5}+4(a_{14}+\epsilon_{4,3})\times (a_{45}+\epsilon_{4,3})}}{2},
	\end{split}
\end{equation*}
where $i=2,3,4$, along with $\tilde{a}_{11}^*=\tilde{a}_{55}^*=1$ and $\tilde{a}_{15}^*=5.4354$, form the best optimally modified PCS. Therefore, $$\tilde{A}_b^*=(1,4.9577,2.3314,1.8975,5.4354), \quad \tilde{A}_w^*=(5.4354,1.0963,2.3314,2.8645,1)^T$$ is the best optimally modified PCS. Thus, by equation \eqref{5equivalence}, $\{0.4270,0.0861,0.1832,0.2251,0.0786\}$ is the best optimal weight set. Now, by equation \eqref{5CR}, we get CR $=0.3119$.\\\\
Table \ref{5example_3} shows that the obtained $\epsilon^*$, interval-weights, and the best optimal weight set coincide with their actual values, validating the proposed framework.\\\\
\textbf{Example 4:} Here, $D_1=\{4\}$, $D_2=\{2,3\}$, and $D_3=\emptyset$. By equations \eqref{5CV} and \eqref{5CV'}, we get $\epsilon_2=0.4355$, $\epsilon_3=0.4401$, $\epsilon_4=0.4865$, $\epsilon_{4,2}=0.0.5480$, and $\epsilon_{4,3}=0.5458$. So, from equation \eqref{5opt_obj_actual}, it follows that $\epsilon^*=\max\{\epsilon_2,\epsilon_3,\epsilon_4,\epsilon_{4,2},\epsilon_{4,3}\}=\epsilon_{4,2}=0.5480$. Now, equation \ref{5fix} gives $\tilde{a}_{15}^*=(a_{14}+\epsilon_{4,2})\times(a_{45}+\epsilon_{4,2})=6.8344$. From equation \eqref{5opt_weights}, we get
\begin{equation*}
	\begin{split}
		{w_1^l}^*&=\frac{\tilde{a}_{15}^*}{1+\tilde{a}_{15}^*+\displaystyle\sum_{j=2,3,4}\min\biggl\{a_{j5}+\epsilon^*,\frac{\tilde{a}_{15}^*}{a_{1j}-\epsilon^*}\biggr\}}\\
		&=\frac{6.8344}{1+6.8344+2.6143+4.0487+2.2801}\\
		&=0.4074.
	\end{split}
\end{equation*}
Similarly, ${w_1^u}^*=0.4077$. Therefore, $w_1^*=[0.4074,0.4077]$. Computing in the same fashion, we obtain $w_2^*=[0.1558,0.1560]$, $w_3^*=[0.2407,0.2413]$, $w_4^*=[0.1359,0.1360]$, and $w_5^*=[0.0596,0.0597]$.\\\\
Since $\epsilon^*=\epsilon_{4,2}$, by statement \textit{3} of Theorem \ref{5best_pcs_actual}, 
\begin{equation*}
	\begin{split}
		&\tilde{a}_{bi}^*=\frac{a_{1i}-a_{i5}+\sqrt{(a_{1i}+a_{i5})^2-4\times a_{1i}\times a_{i5}+4(a_{14}+\epsilon_{4,2})\times (a_{45}+\epsilon_{4,2})}}{2},\\
		&\tilde{a}_{iw}^*=\frac{-a_{1i}+a_{i5}+\sqrt{(a_{1i}+a_{i5})^2-4\times a_{1i}\times a_{i5}+4(a_{14}+\epsilon_{4,2})\times (a_{45}+\epsilon_{4,2})}}{2},
	\end{split}
\end{equation*}
where $i=2,3,4$, along with $\tilde{a}_{11}^*=\tilde{a}_{55}^*=1$ and $\tilde{a}_{15}^*=6.8344$, form the best optimally modified PCS. Therefore, $$\tilde{A}_b^*=(1,2.6143,1.6923,2.9975,6.8344), \quad \tilde{A}_w^*=(6.8344,2.6143,4.0388,2.2801,1)^T$$ is the best optimally modified PCS. Thus, by equation \eqref{5equivalence}, $\{0.4076,0.1559,0.2409,0.1360,0.0596\}$ is the best optimal weight set. Now, by equation \eqref{5CR}, we get CR $=0.1471$.\\\\
Table \ref{5example_4} shows that the obtained $\epsilon^*$, interval-weights, and the best optimal weight set coincide with their actual values, validating the proposed framework.\\\\
\textbf{Example 5:} Here, $D_1=\{2\}$, $D_2=\{3\}$, and $D_3=\emptyset$. By equations \eqref{5CV} and \eqref{5CV'}, we get $\epsilon_2=0.3028$, $\epsilon_3=0.4384$, and $\epsilon_{2,3}=0.5$. So, from equation \eqref{5opt_obj_actual}, it follows that $\epsilon^*=\max\{\epsilon_2,\epsilon_3,\epsilon_{2,3}\}=\epsilon_{2,3}=0.5$. Now, equation \ref{5fix} gives $\tilde{a}_{14}^*=(a_{12}+\epsilon_{2,3})\times(a_{24}+\epsilon_{2,3})=2.25$. From equation \eqref{5opt_weights}, we get
\begin{equation*}
	\begin{split}
		{w_1^l}^*&=\frac{\tilde{a}_{14}^*}{1+\tilde{a}_{14}^*+\displaystyle\sum_{j=2,3}\min\biggl\{a_{j4}+\epsilon^*,\frac{\tilde{a}_{14}^*}{a_{1j}-\epsilon^*}\biggr\}}\\
		&=\frac{2.25}{1+2.25+1.5+1.5}\\
		&=0.36.
	\end{split}
\end{equation*}
Similarly, ${w_1^u}^*=0.36$. Therefore, $w_1^*=[0.36,0.36]$. Computing in the same fashion, we obtain $w_2^*=[0.24,0.24]$, $w_3^*=[0.24,0.24]$, and $w_4^*=[0.16,0.16]$. Therefore, $\{0.36,0.24,0.24,0.16\}$ is the only optimal weight set. Now, by equation \eqref{5CR}, we get CR $=1$.\\\\
Table \ref{5example_5} shows that the obtained $\epsilon^*$ and the optimal weight set coincide with their actual values, validating the proposed framework.\\\\
\textbf{Example 6:} Here, $n_1=2$, $n_2=1$, $D_1=\{3\}$, and $D_2=D_3=\emptyset$. By equation \eqref{5CV}, we get $\epsilon_2=0.1623$. So, from equation \eqref{5opt_obj_multi}, it follows that $\epsilon^*=\epsilon_{2}=0.1623$. Now, equation \ref{5op_PCS_multi} gives $\tilde{a}_{14}^*=\tilde{a}_{24}^*=(a_{13}+\epsilon_2)\times(a_{34}+\epsilon_2)=6.8378$. From equation \eqref{5opt_weights_multi}, we get
\begin{equation*}
	\begin{split}
		{w_1^l}^*={w_2^l}^*&=\frac{\tilde{a}_{14}^*}{1+2\tilde{a}_{14}^*+\min\biggl\{a_{34}+\epsilon^*,\frac{\tilde{a}_{14}^*}{a_{13}-\epsilon^*}\biggr\}}\\
		&=\frac{6.8378}{1+2\times 6.8378+3.1623}\\
		&=0.3833.
	\end{split}
\end{equation*}
Similarly, ${w_1^u}^*={w_2^u}^*=0.3833$. Therefore, $w_1^*=w_2^*=[0.3833,0.3833]$. Computing in the same fashion, we obtain $w_3^*=[0.1773,0.1773]$, and $w_4^*=[0.0561,0.0561]$. Therefore, $\{0.3833,0.3833,0.1773,0.0561\}$ is the only optimal weight set. Now, by equation \eqref{5CR}, we get CR $=0.0436$.\\\\
Note that for $c_1$ and $c_2$, both the interval-weights and the best weights coincide. Additionally, Table \ref{5example_5} indicates that both optimization models \eqref{5optimization_1} and \eqref{5optimization_7} achieve the same optimal objective value. This suggests that optimization model \eqref{5optimization_7} has selected the optimal weight sets in which the weights for $c_1$ and $c_2$ are equal.
\section{Real-World Application}
In this section, we discuss the application of the proposed model in developing a roadmap to address barriers to energy efficiency.\\\\
In the modern era, global energy consumption has increased at an alarming rate. This trend raises significant concerns regarding economic growth, environmental sustainability and the exhaustion of limited energy resources. A renewable energy source capable of fully replacing existing options has yet to be discovered. Therefore, the development and implementation of energy-efficient technologies have become imperative. However, such efforts are often hindered by implementation constraints. These barriers prevent societies from fully leveraging cost-effective and energy-friendly opportunities, making the study of obstacles to energy efficiency a critically important topic. Numerous studies have examined these barriers across various domains including industrial \cite{cagno2015barriers,cagno2013novel}, economic \cite{hirst1990closing,linares2010energy}, and environmental \cite{hanley2009increases,wilkinson2007energy} perspectives among others.\\\\
Buildings represent one of the highest energy-consuming sectors globally. In Europe, for instance, buildings account for a larger share of energy consumption than the industrial sector \cite{zhao2012review,perez2008review}. Energy is consumed throughout a building's entire life cycle including the production of raw materials, construction and demolition. The main reason behind this increasing energy demand and population growth and rising demand for building services \cite{perez2008review}. The actual consumption level varies significantly from one building to another, depending mainly on local weather conditions, structure, and building characteristics \cite{zhao2012review}. Numerous studies have been conducted to predict energy consumption in buildings. Zhao et al. \cite{zhao2012review} comprehensively reviewed various engineering, statistical, and neural network models developed for this purpose. In a similar vein, Caldera et al. \cite{caldera2008energy} analyzed statistical data to evaluate energy demand specifically for space heating in buildings.\\\\
This study aims to rank barriers to energy efficiency in buildings based on their criticality. Gupta et al. \cite{gupta2017developing} identified $27$ such barriers and categorized them into six groups, as presented in Table \ref{5barrier}. The pairwise comparisons between these categories and between barriers within the same category are adopted from their study.\\\\
The category weights, local and global weights of barriers, and the ranking of barriers based on these global weights are presented in Table \ref{5weights_barrier}. Our analysis, based on this model, shows that ``scarcity of financial means'' is the most critical barrier, while ``contention'' is the least critical barrier among all those ranked.
\begin{table}[H]
	\caption {Barriers to energy efficiency in buildings \cite{gupta2017developing}\label{5barrier}}
	\centering
	\begin{tabular}{lll}
		\toprule[0.1em]
		Category &\phantom{a}&Driver\\
		\midrule
		$c_1$: Economic barriers&&$c_{11}$: Scarcity of financial means\\
		&&$c_{12}$: Absence of lucravity\\
		&&$c_{13}$: Poor arbitrage\\
		&&$c_{14}$: Inadequate monetary assessment\\
		&&$c_{15}$: Limits in financial provisioning\\
		$c_2$: Government barriers &&$c_{21}$: Lack of financial motivation \\
		&&$c_{22}$: Bridles in hierarchical inspiration\\
		&&\quad \quad and functional harmony\\
		&&$c_{23}$: Difference in plan of action\\
		&&$c_{24}$: Inappropriate antecedencey\\
		&&$c_{25}$: Lack of standards and references\\
		&&$c_{26}$: Lack of strong authority\\
		$c_3$: Knowledge and learning&&$c_{31}$: Lack of cognizance\\
		\quad \ \ barriers&&$c_{32}$: Inexperient persons\\
		&&$c_{33}$: Lack of information\\
		$c_4$: Market related barriers&&$c_{41}$: Contention\\
		&&$c_{42}$: Perceptual knowledge\\
		&&$c_{43}$: Aspousal\\
		&&$c_{44}$: Dubiety in demand\\
		$c_5$: Organization and social&&$c_{51}$: Cognitive decision making approach\\
		\quad \ \ barriers&&$c_{52}$: Lack of authority and jurisdiction\\
		&&$c_{53}$: Ill defined vision\\
		&&$c_{54}$: Torporpidity in process and practices\\
		$c_6$: Technological barriers&&$c_{61}$: Incompatible technology \\
		&&$c_{62}$: Process related risks\\
		&&$c_{63}$: Lack of energy efficient materials\\
		&&$c_{64}$: No feasibility study\\
		&&$c_{65}$: Slow embodiment of new technology\\
		\bottomrule[0.1em]
	\end{tabular}	
\end{table}
\begin{table}[H]
	\caption {Category weights, barrier weights, and barrier ranking \label{5weights_barrier}}
	\centering
	\begin{tabular}{@{}cccccc@{}}
		\toprule[0.1em]
		Category &Weight&Driver&Local weight&Global weight&Rank\\
		\midrule
		$c_1$&$0.4107$&$c_{11}$&$0.4834$&$0.1985$&$1$\\
		&&$c_{12}$&$0.1196$&$0.0491$&$6$\\
		&&$c_{13}$&$0.0585$&$0.0240$&$13$\\
		&&$c_{14}$&$0.1019$&$0.0419$&$9$\\
		&&$c_{15}$&$0.2366$&$0.0972$&$2$\\
		&&$\epsilon^*$&$0.2583$&-&-\\
		&&CR&$0.0578$&-&-\\
		\midrule
		$c_2$&$0.2225$&$c_{21}$&$0.1440$&$0.0320$&$10$\\
		&&$c_{22}$&$0.2025$&$0.0451$&$8$\\
		&&$c_{23}$&$0.4138$&$0.0920$&$3$\\
		&&$c_{24}$&$0.0501$&$0.0111$&$20$\\
		&&$c_{25}$&$0.0873$&$0.0194$&$17$\\
		&&$c_{26}$&$0.1023$&$0.0228$&$14$\\
		&&$\epsilon^*$&$0.2583$&-&-\\
		&&CR&$0.0578$&-&-\\
		\midrule
		$c_3$&$0.0989$&$c_{31}$&$0.7146$&$0.0707$&$4$\\
		&&$c_{32}$&$0.1$&$0.0099$&$21$\\
		&&$c_{33}$&$0.1854$&$0.0183$&$18$\\
		&&$\epsilon^*$&$0.1459$&-&-\\
		&&CR&$0.0392$&-&-\\
		\midrule
		$c_4$&$0.0459$&$c_{41}$&$0.0713$&$0.0033$&$27$\\
		&&$c_{42}$&$0.5809$&$0.0267$&$12$\\
		&&$c_{43}$&$0.1443$&$0.0066$&$25$\\
		&&$c_{44}$&$0.2035$&$0.0093$&$23$\\
		&&$\epsilon^*$&$0.1459$&-&-\\
		&&CR&$0.0326$&-&-\\
		\midrule
		$c_5$&$0.0847$&$c_{51}$&$0.2605$&$0.0221$&$15$\\
		&&$c_{52}$&$0.1158$&$0.0098$&$22$\\
		&&$c_{53}$&$0.5610$&$0.0475$&$7$\\
		&&$c_{54}$&$0.0627$&$0.0053$&$26$\\
		&&$\epsilon^*$&$0.1538$&-&-\\
		&&CR&$0.0294$&-&-\\
		\midrule
		$c_6$&$0.1373$&$c_{61}$&$0.4724$&$0.0649$&$5$\\
		&&$c_{62}$&$0.1580$&$0.0217$&$16$\\
		&&$c_{63}$&$0.0975$&$0.0134$&$19$\\
		&&$c_{64}$&$0.0528$&$0.0073$&$24$\\
		&&$c_{65}$&$0.2193$&$0.0301$&$11$\\
		&&$\epsilon^*$&$0.1538$&-&-\\
		&&CR&$0.0294$&-&-\\
		\midrule
		$\epsilon^*$&$0.1538$&-&-&-&-\\
		CR&$0.0294$&-&-&-&-\\
		\bottomrule[0.1em]
	\end{tabular}	
\end{table}
\section{Conclusions and Future Directions}
The BWM has emerged as a powerful MCDM technique and is widely adopted in real-world applications. This study introduces an analytical framework for its nonlinear model, which is compatible with any scale and any number of DMs. The framework derives closed-form solutions for optimal interval-weights, CI, and CR. It also identifies the best optimal weight set from the collection of all optimal weight sets. Furthermore, it ranks the barriers to energy efficiency in buildings by their criticality. This research enhances the model's efficiency in several key aspects. By eliminating the need for specialized optimization software, the framework significantly improves computational efficiency. It also rectifies existing inaccuracies in the calculation of CI. Additionally, it transforms CR into an input-based consistency indicator, capable of providing immediate feedback to DMs. The framework also modifies the original optimization model for cases with multiple best or worst criteria to select the optimal weight sets which assign equal weights to all best criteria and equal weights to all worst criteria.\\\\
A promising future research direction emerging from this work is to derive analytical expressions for optimal weights, CI, and CR for models such as the Euclidean BWM \cite{kocak2018euclidean} and the $\alpha$-FBWM \cite{ratandhara2024alpha} for which such closed-form solutions have not yet been established. Developing such frameworks would significantly enhance the applicability and efficiency of these methods.
\section*{Acknowledgements}
The first author gratefully acknowledges the Council of Scientific \& Industrial Research (CSIR), India for financial support to carry out the research work.
\section*{Declaration of Conflict of Interest}
The authors declare that they have no known conflict of financial interests or personal relationships that could have appeared to influence the work reported in this paper.

\end{document}